\documentclass[a4paper, english, 11pt, twoside, reqno]{amsart}

\usepackage{babel}
\usepackage[T1]{fontenc}
\usepackage{amssymb}
\usepackage{amsmath, amsfonts, mathrsfs, amscd}
\usepackage{xcolor}
\usepackage{eucal}
\usepackage[all]{xy}
\usepackage{calligra}
\usepackage[all]{xypic}
\usepackage{graphicx}
\usepackage{hyperref}
\usepackage{amsthm}
\usepackage{pifont}
\usepackage{relsize}

\newtheorem*{theorem-intro}{Theorem}
\newtheorem{theorem}{Theorem}[section]

\newtheorem{proposition}[theorem]{Proposition}
\newtheorem{corollary}[theorem]{Corollary}

\newtheorem{question}{Question}[section]

\theoremstyle{definition}

\newtheorem{remark}[theorem]{Remark}

\def\pf{\begin{proof}}
\def\epf{\end{proof}}

\newcommand{\Na}{\mathbb{N}}
\newcommand{\Z}{\mathbb{Z}}
\newcommand{\Q}{\mathbb{Q}}
\newcommand{\Co}{\mathbb{C}}

\newcommand{\Ow}{\mathcal{O}}

\newcommand\car{\operatorname{char}}

\newcommand{\ot}{\otimes}

\newcommand{\Oint}{\mathcal{O}}

\newcommand{\M}{\mathrm{M}}

\definecolor{azul}{RGB}{0, 47, 103}

\allowdisplaybreaks

\definecolor{rojo}{rgb}{1,0,0}

\DeclareMathAlphabet{\mathpzc}{OT1}{pzc}{m}{it}

\begin{document}

\title[Non-existence of Hopf orders]{Non-existence of Hopf orders for a twist of the alternating and symmetric groups}

\author[J. Cuadra and E. Meir]{Juan Cuadra and Ehud Meir}

\address{J. Cuadra: Universidad de Almer\'\i a, Dpto. Matem\'aticas, 04120 Almer\'\i a, Spain}
\email{jcdiaz@ual.es}

\address{E. Meir: Institute of Mathematics, University of Aberdeen, Fraser Noble Building, Aberdeen AB24 3UE, United Kingdom}
\email{meirehud@gmail.com}

\begin{abstract}
We prove the non-existence of Hopf orders over number rings for two families of complex semisimple Hopf algebras. They are constructed as Drinfel'd twists of group algebras for the following groups: $A_n$, the alternating group on $n$ elements, with $n \geq 5$; and $S_{2m}$, the symmetric group on $2m$ elements, with $m \geq 4$ even. The twist for $A_n$ arises from a $2$-cocycle on the Klein four-group contained in $A_4$. The twist for $S_{2m}$ arises from a $2$-cocycle on a subgroup generated by certain transpositions which is isomorphic to $\Z_2^m$. This provides more examples of complex semisimple Hopf algebras that can not be defined over number rings. As in the previous family known, these Hopf algebras are simple.
\end{abstract}
\maketitle

\section*{Introduction}

The theory of orders lies at the intersection of algebra and number theory. As stated in the preface of Reiner's fundamental book \cite{Re}, ''the beauty of the subject stems from the fascinating interplay between the arithmetical properties of orders and the algebraic properties of the algebras containing them''. \par \smallskip

This interaction is particularly strong in the representation theory of finite groups. For a finite group $G$, the group ring $\Z G$ in the group algebra $\Co G$ is one of the most significant examples of order. An important result here, depending on the existence of such an order, is Frobenius' Theorem: the degree of any complex irreducible representation of $G$ divides the order of $G$, \cite[Proposition 9.32]{CR}. When $\Co G$ is viewed with its customary Hopf algebra structure, $\Z G$ is indeed a Hopf order of $\Co G$. \par \smallskip

Kaplansky predicted that Frobenius' Theorem holds for complex semisimple Hopf algebras, namely: for a complex semisimple Hopf algebra $H$ the dimension of every irreducible representation of $H$ divides the dimension of $H$. This statement, known as Kaplanky's sixth conjecture, remains unanswered. \par \smallskip

Larson introduced Hopf orders in \cite{L} as a first attempt to bring number-theoretic techniques into the representation theory of semisimple Hopf algebras. He proved there that if $H$ admits a Hopf order over a number ring, then $H$ satisfies Kaplansky's conjecture. For a long time it was an open question whether a complex semisimple Hopf algebra always admits a Hopf order over a number ring. We settled this in the negative in \cite{CM} for a family of Hopf algebras studied by Galindo and Natale in \cite{GN}. (Nevertheless, they satisfy the conjecture.) Our result revealed an essential arithmetic difference between group algebras and semisimple Hopf algebras: {\it complex semisimple Hopf algebras may not admit Hopf orders over number rings.} \par \smallskip

In the present paper we further explore this phenomenon. We prove the non-existence of Hopf orders over number rings for two more families of complex semisimple Hopf algebras. Similarly to the previous family, they are constructed as Drinfel'd twists of certain group algebras. The twist is in turn constructed using the following method due to Movshev \cite{Mosh}. \par \smallskip

Let $M$ be an abelian subgroup of $G$. Let $K$ be a number field. Suppose that $K$ is large enough so that the group algebra $K\hspace{-0.9pt}M$ splits. Denote by $\widehat{M}$ the character group of $M$. For $\phi \in  \widehat{M}$ let $e_{\phi}$ be the associated idempotent in $KM$. If $\omega: \widehat{M} \times \widehat{M} \rightarrow K^{\times}$ is a normalized 2-cocycle, then
$$J_{\omega}=\sum_{\phi,\psi \in \widehat{M}} \omega(\phi,\psi) e_{\phi} \otimes e_{\psi}$$
is a twist for $K\hspace{-0.9pt}G$. The twisting procedure alters the coalgebra structure and the antipode of $K\hspace{-0.9pt}G$ but leaves unchanged the algebra structure. So, as algebras, these Hopf algebras are group algebras. \par \smallskip

Here we deal with the following two families of groups and twists:
\begin{enumerate}
\item[(1)] The alternating group $A_n$ on $n$ elements, with $n \geq 5$. Consider the double transpositions $d_1=(12)(34)$ and $d_2=(13)(24)$.\vspace{-1pt} The subgroup $M$ is generated by them. The character group $\widehat{M}$ is generated by $\varphi_i$, for $i=1,2,$ given by $\varphi_i(d_j)=(-1)^{\delta_{ij}}$. The 2-cocycle $\omega$ is the bicharacter:
$$\omega(\varphi_1^{\alpha_1}\varphi_2^{\alpha_2},\varphi_1^{\beta_1}\varphi_2^{\beta_2}) = (-1)^{\alpha_1\beta_2}, \qquad \alpha_i,\beta_j \in \{0,1\}.$$
This family was introduced by Nikshych in \cite{N} and provided the first non-trivial examples of simple and semisimple Hopf algebras. \smallskip

\item[(2)] The symmetric group $S_{2m}$ on $2m$ elements, with $m \geq 4$ even. For $1 \leq i \leq m$ let $t_i$ denote the transposition $(2i-1\hspace{3pt} 2i)$ of $S_{2m}$. The subgroup\vspace{-1pt} $M$ is generated by $t_i$ for $1 \leq i \leq m$. The character group $\widehat{M}$ is generated by $\varphi_i$, for $1 \leq i \leq m,$ given by $\varphi_i(t_j)=(-1)^{\delta_{ij}}$. The 2-cocycle $\omega$ is defined by the formula:
$$\omega(\varphi_1^{\alpha_1}\ldots \varphi_m^{\alpha_m},\varphi_1^{\beta_1}\ldots \varphi_m^{\beta_m}) = (-1)^{\sum_{i<j}\alpha_i\beta_j}, \qquad \alpha_i,\beta_j \in \{0,1\}.$$
This family was introduced by Bichon in \cite{B} and further discussed by Galindo and Natale in \cite{GN}. These examples are also simple.
\end{enumerate}

Our main results, Theorems \ref{main1} and \ref{main2}, are abridged in the following statement:

\begin{theorem-intro}
Let $K$ be a number field with ring of integers $\Oint_K$. Let $R\subset K$ be a Dedekind domain such that $\Oint_K \subseteq R$. Let $G$ be one of the above groups and $J$ the twist arising from the corresponding cocycle. If the twisted Hopf algebra $(K\hspace{-0.9pt}G)_J$ admits a Hopf order over $R$, then $\frac{1}{2}\in R$. \par \smallskip

As a consequence, $(K\hspace{-0.9pt}G)_J$ does not admit a Hopf order over any number ring.
\end{theorem-intro}

This theorem implies that the complexified Hopf algebra $(\Co G)_J$ does not admit a Hopf order over any number ring. \par \smallskip

The strategy of the proof can be outlined as follows. Hopf orders are inherited to Hopf subalgebras and quotient Hopf algebras (Proposition \ref{subsquo}). For our purposes, this fact allows us to focus immediately on $A_5$ in the case of $A_n$. The case of $S_{2m}$ needs several reductions to subgroups and quotient groups to focus ultimately on $S_8$. Assume now that $H$ is a Hopf algebra over $K$ and $X$ is a Hopf order of $H$ over $R$. Consider the dual Hopf algebra $H^*$. The dual order $X^{\star}$ consists of those $\varphi \in H^*$ such that $\varphi(X) \subseteq R$. The dual order is a Hopf order of $H^*$ (Proposition \ref{subsquo}). Another technical result we rely on is that any character of $H$ belongs to $X^{\star}$ and any cocharacter of $H$ belongs to $X$ (Proposition \ref{character}). In particular, any group-like element of $H$ lies in $X$. \par \smallskip \vspace{0.5mm}

As in our previous paper \cite{CM}, the main idea of the proof is to construct elements $\varphi \in X^{\star}$ and $x \in X$, by a subtle manipulation of characters and cocharacters, such that $\varphi(x)=s/2$, with $s$ an odd number. For the characters we directly use the character table of $A_5$ and $S_4$ respectively. For the cocharacters the argument is more involved since the twisting procedure makes the coalgebra structure more difficult to handle. In Proposition \ref{decomp}, we invoke a decomposition of $(K\hspace{-0.9pt}G)_J$ at the coalgebra level by means of the double cosets of the subgroup $M$ in $G$, which was established by Etingof and Gelaki in \cite{EG1}. Given $\tau \in G$ the subspace $K(M \tau M)$ spanned by the double coset $M \tau M$ is a subcoalgebra of $(K\hspace{-0.9pt}G)_J$. We prove in Proposition \ref{simple} that if $M \cap (\tau M \tau^{-1})=\{1\}$, then $K(M\tau M)$ is isomorphic to a matrix coalgebra; and the irreducible cocharacter of $(K\hspace{-0.9pt}G)_J$ associated to it is $|M|e_{\varepsilon}\tau e_{\varepsilon}$. Here $e_{\varepsilon}$ is the idempotent of $K\hspace{-0.9pt}M$ attached to the trivial character. \par \smallskip \vspace{0.5mm}

In the case of $A_5$ the preceding ideas and results do all the work. The case of $S_8$ requires an additional tool, which does not appear in \cite{CM}. We replace the twist $J$ of $K\hspace{-1pt}S_8$ by a cohomologous twist $T$. The Hopf algebras $(K\hspace{-1pt}S_8)_J$ and $(K\hspace{-1pt}S_8)_T$ are isomorphic, but $T$ has a computational advantage:  we are able to show that $T$ and $T^{-1}$ belong to $X \otimes_R X$. Then, $X$ can be twisted by $T^{-1}$ (Proposition \ref{twistorders}), thus obtaining a Hopf order of the group algebra $K\hspace{-1pt}S_8$. Hence $S_8$ is contained in $X$. A permutation of $S_8$ is used to construct the desired element $x$ in $X \cap (K\hspace{-1pt}S_4)$ to which we apply a character of $S_4$ to settle the statement. \par \smallskip \vspace{0.5mm}

These Hopf algebras, as those studied in \cite{CM}, are examples of simple Hopf algebras; i.e., they have no proper normal Hopf subalgebras. In light of our results we wonder if a non-trivial simple and semisimple complex Hopf algebra ever admits a Hopf order over a number ring. In another direction, we ask if a complex semisimple Hopf algebra that is lower-semisolvable, as defined by Montgomery and Witherspoon in \cite{MW}, always admits a Hopf order over a number ring. Some other questions arising from our previous results on orders in semisimple Hopf algebras can be found in \cite[page 2548]{CM} and \cite[page 954]{CM2}. \par \smallskip \vspace{0.5mm}

The paper is organized as follows. Section 1 collects basic material on Hopf orders. The Drinfel'd twist is recalled in Section 2, together with Movshev's method of constructing a twist for a group algebra from a 2-cocycle on an abelian subgroup. Here we present the coalgebra decomposition of the twisted group algebra in terms of double cosets. We also describe the irreducible cocharacter for a simple subcoalgebra provided by a double coset. In Section 3 we prove the non-existence of Hopf orders for the twist of the alternating group. For the symmetric group this is done in Section 4. Several final questions are formulated in Section 5. \newpage

\section{Preliminaries}\label{prelim}
\setcounter{equation}{0}

In this section we fix notation and gather together several results on Hopf orders from \cite[Subsection 1.2]{CM} that we will use later. We refer the reader to there for their proofs. \par \smallskip

\subsection{Notation} Throughout $H$ is a finite-dimensional Hopf algebra over a base field $K$. Vector spaces, linear maps, and unadorned tensor products are over $K$, unless otherwise indicated. The identity element of $H$ is denoted by $1_H$ and the coproduct, counit, and antipode by $\Delta, \varepsilon,$ and $\mathcal{S}$ respectively. The dual Hopf algebra of $H$ is denoted by $H^*$. For $\varphi \in H^*$ and $h \in H$ we sometimes use the duality pairing notation $\langle \varphi, h \rangle$ instead of $\varphi(h)$. Our references for the theory of Hopf algebras are \cite{Mo} and \cite{Ra}. \par \smallskip

\subsection{Hopf orders} Let $R$ be a subring of $K$ and $V$ a finite-dimensional vector space over $K$. An \emph{order of\hspace{2pt} $V$ over $R$} is a finitely generated and projective $R$-submodule $X$ of $V$ such that the natural map $X \otimes_R K \rightarrow V$ is an isomorphism. The submodule $X$ corresponds to the image of $X \otimes_R R$. \par \smallskip

A \textit{Hopf order of $H$ over $R$} is an order $X$ of $H$ such that $1_H \in X$, $XX \subseteq X$, $\Delta(X)\subseteq X\otimes_{R} X$, $\varepsilon(X) \subseteq R,$ and $\mathcal{S}(X)\subseteq X$. (For the coproduct, $X\otimes_{R} X$ is naturally identified as an $R$-submodule of $H\otimes H$.) Then, $X$ is a Hopf algebra over $R$, which is finitely generated and projective as an $R$-module, such that $X\otimes_{R} K \simeq H$ as Hopf algebras over $K$. \par \smallskip

In the next two results $K$ is a number field and $R\subset K$ is a Dedekind domain containing the ring of algebraic integers $\Oint_K$. Hopf orders are over $R$. \par \smallskip

\begin{proposition}\cite[Propositions 1.1 and 1.9]{CM} \label{subsquo} Let $X$ be a Hopf order of $H$. \vspace{2pt}
\begin{enumerate}
\item[{\it (i)}] The dual order $X^{\star}:=\{\varphi \in H^* : \varphi(X) \subseteq R\}$ is a Hopf order of $H^*$. \vspace{2pt}
\item[{\it (ii)}] The natural isomorphism $H \simeq H^{**}$ induces an isomorphism of Hopf orders $X \simeq X^{\star \star}$. \vspace{2pt}
\item[{\it (iii)}] If $A$ is a Hopf subalgebra of $H$, then $X\cap A$ is a Hopf order of $A$. \vspace{2pt}
\item[{\it (iv)}] If $\pi:H \rightarrow B$ is a surjective Hopf algebra map, then $\pi(X)$ is a Hopf order of $B$.
\end{enumerate}
\end{proposition}

The proofs of our main results are based on the following:

\begin{proposition}\cite[Proposition 1.2]{CM} \label{character}
Let $X$ be a Hopf order of $H$. Any character of $H$ belongs to $X^{\star}$ and any character of $H^*$ (cocharacter of $H$) belongs to $X$. In particular, any group-like element of $H$ belongs to $X$.
\end{proposition}

\section{Coalgebra structure of twistings arising from abelian subgroups}
\setcounter{equation}{0}

We discuss here Movshev's idea \cite[Section 1]{Mosh} of deforming a group algebra through a $2$-cocycle on an abelian subgroup and the link between the deformed coalgebra structure and the double cosets of the subgroup set up by Etingof and Gelaki in \cite[Section 3]{EG1}. \par \smallskip

We start by recalling the general deformation procedure, due to Drinfel'd. \vspace{1pt} An invertible element $J:=\sum J^{(1)} \otimes J^{(2)} \in H \otimes H$ is a {\it twist} for $H$ provided that:
$$\begin{array}{c}
(1_H \otimes J)(Id \otimes \Delta)(J)=(J \otimes 1_H)(\Delta \otimes Id)(J), \quad \textrm{and} \vspace{4pt} \\
(\varepsilon \otimes Id)(J)=(Id  \otimes \varepsilon)(J)=1_H.
\end{array}$$
The {\it Drinfel'd twist} of $H$ is the new Hopf algebra $H_J$ defined as follows: $H_J=H$ as an algebra, the counit is that of $H$, and the new coproduct and antipode are given by:
$$\Delta_J(h)=J \Delta(h)J^{-1} \qquad \textrm{and} \qquad \mathcal{S}_J(h)=U_J\hspace{1pt}\mathcal{S}(h)\hspace{0.5pt}U_J^{-1} \qquad \forall h \in H.$$
Here $U_J:=\sum J^{(1)}\mathcal{S}(J^{(2)})$ and $U_J^{-1}=\sum \mathcal{S}(J^{-(1)})J^{-(2)},$ \vspace{1pt} where for the latter we write $J^{-1}=\sum J^{-(1)} \otimes J^{-(2)}$. \par \smallskip

Two twists $J$ and $J'$ for $H$ are said to be {\it cohomologous} if there is $v \in H$ invertible such that $J'=(v \otimes v)J\Delta(v^{-1})$. In this case, the map $f:H_{J} \rightarrow H_{J'}, \, h \mapsto vhv^{-1}$ is a Hopf algebra isomorphism. \par \smallskip

Let $G$ be a finite group and consider the group algebra $K\hspace{-0.9pt}G$. \vspace{1pt} Let $M$ be an abelian subgroup of $G$. Suppose that $\car K \nmid \vert M \vert$ and that $K$ is large enough so that $K\hspace{-0.9pt}M$ splits. Denote by $\widehat{M}$ the character group of $M$. For $\phi \in  \widehat{M}$, the corresponding idempotent in $K\hspace{-1pt}M$ is given by:
$$e_{\phi} = \frac{1}{\vert M \vert } \sum_{m \in M} \phi(m^{-1}) m.$$
If $\omega: \widehat{M} \times \widehat{M} \rightarrow K^{\times}$ is a normalized 2-cocycle, then
$$J_{\omega}=\sum_{\phi,\psi \in \widehat{M}} \omega(\phi,\psi) e_{\phi} \otimes e_{\psi}$$
is a twist for $K\hspace{-0.9pt}M$. Since $K\hspace{-0.9pt}M$ is a Hopf subalgebra of $K\hspace{-0.9pt}G$, it is also a twist for $K\hspace{-0.9pt}G$. This way of twisting a group algebra plays a key role in the classification of triangular and cotriangular Hopf algebras (\cite{EG2}, \cite{AEGN}, and \cite{EOV}) and in the construction of simple and semisimple Hopf algebras (\cite{N} and \cite{GN}). \par \smallskip

If $\omega$ and $\omega'$ are cohomologous cocycles, then the corresponding twists $J_{\omega}$ and $J_{\omega'}$ are cohomologous in the above sense. Concretely, let $q:\widehat{M} \rightarrow K^{\times}$ be a map such that $\omega'=\omega \partial(q)$, where $\partial$ is the coboundary map. The element $v$ is given by $v=\sum_{\phi \in \widehat{M}} q(\phi)e_{\phi}.$ \par \smallskip

In the sequel we just write $J$ for the twist $J_{\omega}$. The following result delves into the coalgebra structure of $(K\hspace{-0.9pt}G)_J$. The first item is \cite[Proposition 3.1]{EG1} and the second one is a reinterpretation of \cite[Proposition 4.1]{EG1}. The proofs are included since they contain constructions and formulas that we will use later.

\begin{proposition}\label{decomp}
Let $\{\tau_{\ell}\}_{\ell \in \Lambda}$ be a set of representatives of the double cosets of $M$ in $G$. Then:
\begin{enumerate}
\item[{\it (i)}] As a coalgebra, $(K\hspace{-0.9pt}G)_J$ decomposes as the direct sum of subcoalgebras
\begin{equation}\label{decompKGJ}
(K\hspace{-0.9pt}G)_J = \bigoplus_{\ell \in \Lambda} K(M\tau_{\ell} M).
\end{equation}

\item[{\it (ii)}] The dual algebra of $K(M\tau_{\ell} M)$ is isomorphic to the twisted group algebra $K^{(\omega,\omega^{-1})\vert _{N_{\ell}}}[N_{\ell}]$, where
$$\hspace{1cm} N_{\ell}=\big\{(\phi,\psi) \in \widehat{M}\times\widehat{M} \ : \ \psi(m)=\phi(\tau_{\ell} m \tau_{\ell}^{-1}) \hspace{7pt} \forall m \in M\cap (\tau_{\ell} M \tau_{\ell}^{-1}) \big\}.$$
In particular, if $M \cap (\tau_{\ell} M \tau_{\ell}^{-1}) =\{1\}$, then $N_{\ell}=\widehat{M} \times\widehat{M}$.
\end{enumerate}
\end{proposition}

\pf (i) That \eqref{decompKGJ} is the direct sum of those subspaces follows from the fact that $\{M\tau_{\ell} M\}_{\ell \in \Lambda}$ is a partition of $G$. To see that $K(M\tau_{\ell} M)$ is a subcoalgebra of $(K\hspace{-0.9pt}G)_J$ bear in mind that $\Delta_J(\sigma) = J(\sigma \otimes \sigma)J^{-1}$ for every $\sigma \in M\tau_{\ell} M$ and $J$ is supported on $M \otimes M$. Nevertheless, we will compute explicitly the coproduct of $e_{\phi}\tau_{\ell} e_{\psi}$ for a later application:
\begin{align}
\Delta_J(e_{\phi}\tau_{\ell} e_{\psi})  & = J \Delta(e_{\phi}\tau_{\ell} e_{\psi})J^{-1} \notag \vspace{5pt} \\
               & \overset{\text{\ding{173}}}{=} {\displaystyle J \bigg(\sum_{\lambda,\rho \in \widehat{M}} e_{\lambda} \tau_{\ell} e_{\rho} \otimes e_{\lambda^{-1}\phi} \tau_{\ell} e_{\rho^{-1}\psi}\bigg)J^{-1}} \notag \vspace{5pt} \\
               & = {\displaystyle \sum_{\lambda,\rho \in \widehat{M}} \hspace{5pt}  \sum_{\phi',\psi' \in \widehat{M}} \hspace{5pt} \sum_{\phi'',\psi'' \in \widehat{M}}  \omega(\phi',\psi') \omega^{-1}(\phi'',\psi'') e_{\phi'} e_{\lambda} \tau_{\ell} e_{\rho} e_{\phi''}  } \notag \vspace{-12pt} \\
               &  \hspace{7.8cm} {\displaystyle \otimes e_{\psi'} e_{\lambda^{-1}\phi} \tau_{\ell} e_{\rho^{-1}\psi} e_{\psi''}} \notag \vspace{8pt} \\
               & \overset{\text{\ding{172}}}{=} {\displaystyle \sum_{\lambda,\rho \in \widehat{M}} \omega(\lambda,\lambda^{-1}\phi) \omega^{-1}(\rho,\rho^{-1}\psi) e_{\lambda} \tau_{\ell} e_{\rho} \otimes e_{\lambda^{-1}\phi} \tau_{\ell} e_{\rho^{-1}\psi}.} \label{eqcopr}
\end{align}
We used here that:
\begin{enumerate}
\item[\ding{172}] $\{e_{\phi}\}_{\phi \in \widehat{M}}$ is a complete set of orthogonal idempotents in $K\hspace{-0.9pt}M$; \vspace{3pt}
\item[\ding{173}] $\Delta(e_{\phi})=\sum\limits_{\lambda \in \widehat{M}} e_{\lambda} \otimes e_{\lambda^{-1}\phi}$.
\end{enumerate}
\par \medskip

(ii) The dimension of $K(M\tau_{\ell} M)$ equals the cardinal of the double coset $M\tau_{\ell} M$. This is given by the formula:
$$\vert M\tau_{\ell} M \vert = \frac{\vert M \vert^2}{\vert M \cap (\tau_{\ell} M \tau_{\ell}^{-1})\vert}.$$
Notice that $K(M\tau_{\ell}M)$ is spanned by $\{e_{\phi}\tau_{\ell}e_{\psi}\}_{(\phi,\psi) \in \widehat{M} \times \widehat{M}}$. For
$m \in M \cap (\tau_{\ell} M \tau_{\ell}^{-1})$ it holds that
$$\psi(m)e_{\phi}\tau_{\ell}e_{\psi} = e_{\phi}\tau_{\ell}me_{\psi}= \phi(\tau_{\ell}m\tau_{\ell}^{-1})e_{\phi}\tau_{\ell}e_{\psi}.$$
If $(\phi,\psi) \notin N_{\ell}$, there is $m\in M \cap (\tau_{\ell} M \tau_{\ell}^{-1})$ such that $\psi(m) \neq \phi(\tau_{\ell}m\tau_{\ell}^{-1})$.\vspace{2pt} Then $e_{\phi}\tau_{\ell}e_{\psi}=0$. Therefore, $\{e_{\phi}\tau_{\ell}e_{\psi}\}_{(\phi,\psi)\in N_{\ell}}$ spans $K(M\tau_{\ell}M)$. On the other hand, identify naturally $\widehat{M}\times\widehat{M}$ with $\widehat{M \times M}$. Consider the following subgroup $L$ of $M\times M$:
$$L=\big\{(m,\tau_{\ell}m\tau_{\ell}^{-1})\ : \ m\in M \cap (\tau_{\ell} M \tau_{\ell}^{-1})\big\}.$$
Then, $(\phi,\psi) \in N_{\ell}$ if and only if $(\psi,\phi^{-1}) \in L^{\perp}$, the perpendicular subgroup of $L$. We have:
$$\vert N_{\ell} \vert = \vert L^{\perp} \vert = \frac{\vert M \times M\vert}{\vert L \vert}=\frac{|M|^2}{\vert M \cap (\tau_{\ell} M \tau_{\ell}^{-1}) \vert }= \vert M\tau_{\ell}M \vert.\vspace{2pt}$$
This implies that $\{e_{\phi}\tau_{\ell}e_{\psi}\}_{(\phi,\psi)\in N_{\ell}}$ is a basis of $K(M\tau_{\ell}M)$. \par \smallskip

Let $\{u_{(\phi,\psi)}\}_{(\phi,\psi)\in N_{\ell}}$ denote the dual basis of $\{e_{\phi}\tau_{\ell}e_{\psi}\}_{(\phi,\psi)\in N_{\ell}}$ in $K(M\tau_{\ell} M)^*$. We compute the product of two such elements:
$$\begin{array}{l}
\big\langle u_{(\phi_1,\psi_1)} u_{(\phi_2,\psi_2)}, e_{\phi}\tau_{\ell} e_{\psi} \big\rangle = \big\langle u_{(\phi_1,\psi_1)} \otimes u_{(\phi_2,\psi_2)}, \Delta_J(e_{\phi}\tau_{\ell} e_{\psi}) \big\rangle \vspace{5pt} \\
\hspace{1.24cm} {\displaystyle \overset{\text{\eqref{eqcopr}}}{=} \sum_{\lambda,\rho \in \widehat{M}} \omega(\lambda,\lambda^{-1}\phi) \omega^{-1}(\rho,\rho^{-1}\psi) \big\langle u_{(\phi_1,\psi_1)}, e_{\lambda} \tau_{\ell} e_{\rho} \big\rangle \big\langle u_{(\phi_2,\psi_2)}, e_{\lambda^{-1}\phi} \tau_{\ell} e_{\rho^{-1}\psi} \big\rangle }\vspace{5pt} \\
\hspace{1.4cm} {\displaystyle = \omega(\phi_1,\phi_2) \omega^{-1}(\psi_1,\psi_2) \delta_{\phi_1\phi_2,\phi} \delta_{\psi_1\psi_2,\psi}} \vspace{5pt} \\
\hspace{1.4cm} {\displaystyle = \big\langle \omega(\phi_1,\phi_2) \omega^{-1}(\psi_1,\psi_2)u_{(\phi_1\phi_2,\psi_1\psi_2)}, e_{\phi}\tau_{\ell} e_{\psi} \big\rangle . }
\end{array}$$
Hence,
$$u_{(\phi_1,\psi_1)} u_{(\phi_2,\psi_2)}=\omega(\phi_1,\phi_2) \omega^{-1}(\psi_1,\psi_2)u_{(\phi_1\phi_2,\psi_1\psi_2)}.$$
This is precisely the product of $K^{(\omega,\omega^{-1})\vert _{N_{\ell}}}[N_{\ell}]$.
\epf

We will need the following known facts on the twisted group algebra: \par \smallskip

Let $N$ be a finite abelian group. \hspace{-1.2pt}Recall that a normalized $2$-cocycle \mbox{$c\hspace{-0.45pt}: \hspace{-0.45pt} N \hspace{-2pt}\times \hspace{-1pt} N \hspace{-0.45pt} \rightarrow \hspace{-0.45pt} K^{\times}$} is \textit{non-degenerate} if the skew-symmetric bilinear\vspace{-0.25pt} form $c\hspace{1pt} c_{21}^{-1}\hspace{-0.5pt}: \hspace{-0.5pt} N \hspace{-1.5pt}\times \hspace{-0.5pt} N \hspace{-0.25pt} \rightarrow \hspace{-0.25pt} K^{\times},$ \linebreak $(g,h) \mapsto c(g,h)c^{-1}(h,g)$ \hspace{1pt}is \hspace{1pt}non-degenerate. An \hspace{1pt}element $g \in N$ is \textit{$c$-regular} if \linebreak $c(g,h)=c(h,g)$ for all $h \in N$. The set consisting of $c$-regular elements, denoted by $Rad(c)$, is a subgroup of $N$. \par \smallskip

Suppose that $\car K \nmid \vert N \vert$. The twisted group algebra $K^cN=\oplus_{g \in N} K\hspace{-1pt} u_g$ is then semisimple. Assume that $K$ is large enough so that $K^c N$ splits as an algebra. The center of $K^c N$ is spanned, as a vector space, by the set $\{u_g: g \ \textrm{is}\  c\textrm{-regular}\}$. Thus, the following statements are equivalent:
\begin{enumerate}
\item[(1)] $K^c N$ is simple;
\item[(2)] $1$ is the only $c$-regular element;
\item[(3)] $c$ is non-degenerate.
\end{enumerate}

The algebra $K^c N$ decomposes as a direct sum of $|Rad(c)|$ matrix algebras, each of dimension $|N/Rad(c)|$. Hence, all irreducible representations of $K^c N$ are of dimension $\sqrt{|N/Rad(c)|}$. The details about this can be found, for example, in \cite[Section 9.1]{MM}. \par \smallskip

The Hopf algebra $(K\hspace{-0.9pt}G)_J$ is cosemisimple by \cite[Corollary 3.6]{AEGN}. \vspace{1pt} There is a special case in which $K(M\tau_{\ell} M)$ provides a simple component of its Wedderburn decomposition. We assume that $K$ is sufficiently large so that $(K\hspace{-0.9pt}G)_J$ splits as a direct sum of matrix coalgebras.

\begin{proposition}\label{simple}
Let $M$ and $\tau_{\ell}$ be as before. If $M \cap (\tau_{\ell} M \tau_{\ell}^{-1}) =\{1\}$ and $\omega$ is non-degenerate, then $K(M\tau_{\ell} M)$ is isomorphic to a matrix coalgebra of size $\vert M \vert.$ Moreover, the irreducible cocharacter of $(K\hspace{-0.9pt}G)_J$ attached to $K(M\tau_{\ell} M)$ is $|M|e_{\varepsilon}\tau_{\ell}e_{\varepsilon}$.
\end{proposition}

\pf Since $M \cap (\tau_{\ell} M \tau_{\ell}^{-1}) =\{1\}$, by the above result, $K(M\tau_{\ell} M)^*$\vspace{0.5pt} is isomorphic to $K^{(\omega,\omega^{-1})}[\widehat{M} \times \widehat{M}\hspace{2pt}]$. This is in turn isomorphic to $K^{\omega}[\widehat{M}\hspace{1pt}] \otimes K^{\omega^{-1}}[\widehat{M}\hspace{1pt}]$.\vspace{1pt} Under the assumption on $K$ and $\omega$, the algebra $K^{\omega}[\widehat{M}\hspace{1pt}]$\vspace{1pt} is isomorphic to a matrix algebra. Hence $K^{\omega}[\widehat{M}\hspace{1pt}] \otimes K^{\omega^{-1}}[\widehat{M}\hspace{1pt}]$ is so as well. \par \smallskip

We next prove the second statement. The character of the regular representation of $K^{(\omega,\omega^{-1})}[\widehat{M} \times \widehat{M}\hspace{2pt}]$ maps
$$u_{(\phi,\psi)} \mapsto
\begin{cases}
|M|^2 & \textrm{if }\  (\phi,\psi)=(\varepsilon,\varepsilon), \\
\hspace{7pt} 0   & \textrm{if }\  (\phi,\psi)\neq (\varepsilon,\varepsilon).
\end{cases}$$
This is because for $(\phi',\psi') \neq (\varepsilon,\varepsilon)$, the multiplication by $u_{(\phi',\psi')}$ permutes (up to a non-zero scalar) the basis $\{u_{(\phi,\psi)}\}_{\phi, \psi \in \widehat{M}}$ without fixed points, and $|M|^2$ is the dimension of the algebra.
The regular representation of a matrix algebra of size $\vert M \vert$ is the direct sum of $\vert M \vert$ copies of the unique irreducible representation. Then, the character $\Gamma$ of the irreducible representation of $K^{(\omega,\omega^{-1})}[\widehat{M} \times \widehat{M}\hspace{2pt}]$ is given by
$$\Gamma(u_{(\phi,\psi)})= \vert M \vert \delta_{\phi,\varepsilon}\delta_{\psi,\varepsilon}, \qquad \forall \phi, \psi \in \widehat{M}.$$
Recall that $\{u_{(\phi,\psi)}\}_{\phi, \psi \in \widehat{M}}$ is the dual basis of $\{e_{\phi}\tau_{\ell} e_{\psi}\}_{\phi,\psi \in \widehat{M}}$.\vspace{1pt}  Identifying $K(M\tau_{\ell} M)^{**}$ and $K(M\tau_{\ell} M)$ via these bases, $\Gamma$ corresponds to $\vert M \vert e_{\varepsilon}\tau_{\ell} e_{\varepsilon}$.
\epf

\begin{remark}
Our previous discussion on the corepresentations of $(K\hspace{-0.9pt}G)_J$ only includes those aspects needed for our purposes. The irreducible corepresentations of any subcoalgebra $K(M\tau_{\ell} M)$ are described by Etingof and Gelaki in \cite[Theorem 3.2]{EG1}. The proof of that result requires deeper arguments. We refer the interested reader to there.
\end{remark}

The last result of this part relates the twisting procedure and Hopf orders. It will be crucial in establishing Theorem \ref{main2}.

\begin{proposition}\label{twistorders}
Let $H$ be a Hopf algebra over $K$ and $J$ a twist for $H$. Let $R$ be a subring of $K$ and $X$ a Hopf order of $H$ over $R$. Assume that $J$ and $J^{-1}$ belong to $X \otimes_R X$. Then, $X$ is a Hopf order of $H_J$.
\end{proposition}

\pf We must show that $X$ is closed under the operations of $H_J$. Since the algebra structure of $H_J$ remains unchanged, we only need to deal with the coalgebra structure and the antipode. For the former, notice that the hypothesis $J^{\pm 1} \in X \otimes_R X$ implies $\Delta_J(X) = J\Delta(X)J^{-1} \subset X\otimes_R X$. For the latter, notice that $U_J, U_J^{-1} \in X$ and $\mathcal{S}_J(X) = U_J\hspace{1pt}\mathcal{S}(X)\hspace{0.5pt}U_J^{-1} \subseteq X$.
\epf

As done for Hopf algebras, we write $X_J$ for the twisted Hopf order.

\section{Twisting of the alternating group}
\setcounter{equation}{0}

We next illustrate Propositions \ref{decomp} and \ref{simple} with the twisting of the alternating group studied by Nikshych in \cite[Subsection 5.2]{N}. \par \smallskip

For $n \geq 4$ let $A_n$ be the alternating group on the set $\{1,\ldots,n\}$. Consider the abelian subgroup $M=\{id,(12)(34),(13)(24),(14)(23)\}$.\vspace{-1pt}  Choose $a=(12)(34)$ and $b=(13)(24)$ as generators.\vspace{-1pt} The character group $\widehat{M}$ consists of the following elements:

$$\begin{array}{c|rrrr}
             &  1  &  a  &  b  & ab \\
\hline
\varepsilon  &  1  &   1  &   1  & 1   \\
\varphi_1    &  1  &   -1  &  1  & - 1 \\
\varphi_2    &  1  &  1   &   -1  & - 1 \\
\varphi_1 \varphi_2  &  1  &  -1  &  -1  &   1
\end{array}$$
\smallskip

The orthogonal idempotents giving the Wedderburn decomposition of $K\hspace{-0.9pt}M$ are:
$$\begin{array}{lp{1cm}l}
{\displaystyle e_{\varepsilon} = \frac{1}{4}(id+a+b+ab),} &  & {\displaystyle e_{\varphi_1} = \frac{1}{4}(id-a+b-ab),} \vspace{5pt} \\
{\displaystyle e_{\varphi_2} = \frac{1}{4}(id+a-b-ab),}        &  & {\displaystyle e_{\varphi_1\varphi_2} = \frac{1}{4}(id-a-b+ab).} \\
\end{array}$$

Assume that $K$ contains a primitive 4th root of unity $\xi$ and consider the cocycle \linebreak $\omega: \widehat{M} \times \widehat{M} \rightarrow K^{\times}$ defined by the following table:

$$\begin{array}{c|rrrr}
             &  \varepsilon & \varphi_1 & \varphi_2 & \varphi_1\varphi_2 \\
\hline
\varepsilon        &  1  &  1   &  1     &  1   \\
\varphi_1          &  1  &  1   &  \xi   & -\xi  \\
\varphi_2          &  1  & -\xi  &  1     & \xi  \\
\varphi_1\varphi_2 &  1  & \xi   & -\xi   &  1
\end{array}$$
\smallskip

Take the twist
\begin{equation}\label{cocAn}
J=\sum_{\phi,\psi \in \widehat{M}} \omega(\phi,\psi) e_{\phi} \otimes e_{\psi}
\end{equation}
for $K\hspace{-0.9pt}M$ afforded by $\omega$.
\par \smallskip

We compute the Wedderburn decomposition, at the coalgebra level, of the twisted Hopf algebra $(K\hspace{-1pt}A_5)_J.$

\begin{proposition}\label{A5Jdecomp}
Let $\mathcal{T}$ denote the subset $\{(12345),(13524),(23)(45)\}$ of $A_5$. Then:
\begin{enumerate}
\item[{\it (i)}] The Hopf algebra $(K\hspace{-1pt}A_5)_J$ decomposes as a direct sum of simple subcoalgebras:
\[(K\hspace{-1pt}A_5)_J=\Big(\bigoplus_{g \in A_4} Kg\Big) \scalebox{1.1}{$\bigoplus$} \Big(\bigoplus_{\tau \in \mathcal{T}} K(M\tau M)\Big).\]
All subcoalgebras $K(M\tau M)$ are isomorphic to the matrix coalgebra $\M_4^c(K)$. \vspace{3pt}

\item[{\it (ii)}] For each $\tau \in \mathcal{T}$, the irreducible cocharacter of $(K\hspace{-1pt}A_5)_J$ arising from $K(M\tau M)$ is:
$$\mu_{\tau}=4e_{\varepsilon}\tau e_{\varepsilon}.$$

\item[{\it (iii)}] The coproduct of $\mu_{\tau}$ is:
$$\Delta_J(\mu_{\tau}) = \Delta(\mu_{\tau})= \hspace{5pt}{\displaystyle \frac{1}{4}\sum_{\sigma \in M \tau M} \sigma \otimes \sigma.}$$
\end{enumerate}
\end{proposition}

\pf (i) We particularize the decomposition \eqref{decompKGJ} to this example. We must determine all double cosets of $M$ in $A_5$. We start by analyzing the case of an element $g \in A_4$. Since $M$ is normal in $A_4$, the double coset $MgM$ is just the coset $Mg$. We next verify that $Kh$ is a subcoalgebra of $(K\hspace{-1pt}A_5)_J$ for every $h \in Mg$:
$$\begin{array}{l}
\hspace{-1mm}\Delta_J(h)  = J \Delta(h)J^{-1} \vspace{5pt} \\
\hspace{4mm} = {\displaystyle \sum_{\phi,\psi \in \widehat{M}} \hspace{5pt} \sum_{\phi',\psi' \in \widehat{M}} \omega(\phi,\psi) \omega^{-1}(\phi',\psi') h (h^{-1}e_{\phi} h) e_{\phi'} \otimes h (h^{-1} e_{\psi} h) e_{\psi'}} \vspace{5pt} \\
\hspace{4mm} \overset{\text{\ding{172}}}{=} {\displaystyle \sum_{\phi',\psi' \in \widehat{M}} \hspace{5pt} \sum_{\phi'',\psi'' \in \widehat{M}} \omega(h \succ \phi'' \prec h^{-1}, h \succ \psi'' \prec h^{-1}) \omega^{-1}(\phi',\psi') h e_{\phi''} e_{\phi'} \otimes h e_{\psi''} e_{\psi'}} \vspace{5pt} \\
\hspace{4mm} \overset{\text{\ding{173}}}{=} {\displaystyle \sum_{\phi',\psi' \in \widehat{M}} \hspace{5pt} \sum_{\phi'',\psi'' \in \widehat{M}} \omega(\phi'', \psi'') \omega^{-1}(\phi',\psi') h e_{\phi''} e_{\phi'} \otimes h e_{\psi''} e_{\psi'}} \vspace{5pt} \\
\hspace{4mm} = {\displaystyle \sum_{\phi,\psi \in \widehat{M}} \omega(\phi, \psi) \omega^{-1}(\phi,\psi) h e_{\phi} \otimes h e_{\psi}} \vspace{5pt} \\
\hspace{4mm} = h \otimes h.
\end{array}$$
In this computation we have used:
\begin{enumerate}
\item[\ding{172}] The idempotent $h^{-1}e_{\phi} h$ must be some $e_{\phi''}$ for a unique $\phi'' \in \widehat{M}$. Indeed, $\phi''=h^{-1} \succ \phi \prec h$, where $\langle h^{-1} \succ \phi \prec h , m \rangle = \langle \phi, hmh^{-1} \rangle$ for all $m \in M;$ \vspace{3pt}
\item[\ding{173}] The cocycle $\omega$ is invariant under conjugation by elements in $A_4$; namely,\vspace{-5pt}
$$\omega(h \succ \phi'' \prec h^{-1}, h \succ \psi'' \prec h^{-1})=\omega(\phi'', \psi''), \qquad \forall \phi'', \psi'' \in \widehat{M}.$$
\end{enumerate}
(This shows in fact that $K\hspace{-1pt}A_4$ remains unchanged when twisted by $J$.)
\par \smallskip

On the other hand, notice that $M\cap (tMt^{-1})=\{id\}$ for every $t \in A_5 \backslash A_4$. Hence, each double coset outside $A_4$ contains exactly $16$ elements. This implies that there are $3$ double cosets outside $A_4$. We can choose as a set of representatives $$\mathcal{T}=\{(12345),(13524),(23)(45)\}.$$

Since $\omega$ is non-degenerate, by Proposition \ref{simple}, $K(M\tau M)$ is isomorphic to the matrix coalgebra $\M_4^c(K)$. \par \medskip

(ii) Apply also Proposition \ref{simple}. \par \medskip

(iii) Finally, we compute $\Delta_J(\mu_{\tau})$:
$$\begin{array}{rl}
\hspace{1cm} \Delta_J(\mu_{\tau})  & \hspace{-2mm} = 4 \Delta_J(e_{\varepsilon}\tau e_{\varepsilon}) \vspace{5pt} \\
\phantom{\Delta_J(\mu_{\tau})}  & \hspace{-7.5mm} \overset{\text{\eqref{eqcopr}}}{=} {\displaystyle  4 \sum_{\lambda,\rho \in \widehat{M}} \omega(\lambda,\lambda^{-1})\omega^{-1}(\rho,\rho^{-1})e_{\lambda}\tau e_{\rho}\otimes e_{\lambda^{-1}}\tau e_{\rho^{-1}}} \vspace{5pt} \\
  & \hspace{-6mm} = {\displaystyle 4\sum_{\lambda,\rho\in \widehat{M}} e_{\lambda}\tau e_{\rho}\otimes e_{\lambda}\tau e_{\rho}}  \hspace{1.1cm} \big(\textrm{as $\lambda=\lambda^{-1}$ and $\omega(\lambda,\lambda)=1$ \hspace{1.5pt} $\forall \lambda\in\widehat{M}$}\hspace{1.5pt}\big) \vspace{5pt} \\
  & \hspace{-6mm} = \Delta(4 e_{\varepsilon}\tau e_{\varepsilon}) \vspace{5pt} \\
  & \hspace{-6mm} = \Delta (\mu_{\tau}).
\end{array}\hspace{-2mm}$$

Substituting $e_{\varepsilon}= \frac{1}{4}\sum_{k,l}a^kb^l$ in $\mu_{\tau}=4e_{\varepsilon}\tau e_{\varepsilon}$ we get the expression
$$\mu_{\tau}= \frac{1}{4}\sum_{k,l,r,s} a^kb^l \tau a^rb^s.$$
Using the previous computation, we arrive at:
$$\begin{array}{rl}
\Delta_J(\mu_{\tau}) & \hspace{-2mm} = \Delta (\mu_{\tau}) \vspace{5pt} \\
  & \hspace{-2mm} = {\displaystyle \frac{1}{4} \sum_{k,l,r,s} a^kb^l \tau a^rb^s \otimes a^kb^l \tau a^rb^s} \vspace{5pt} \\
  & \hspace{-2mm} = {\displaystyle \frac{1}{4}\sum_{\sigma \in M \tau M} \sigma \otimes \sigma.}\vspace{-5pt}
\end{array}$$
\epf

\begin{remark}
The previous proof uses in an essential way the concrete form of the cocycle $\omega$. This was our reason for doing this choice. However, for our purposes, we do not lose generality since cohomologous cocycles produce isomorphic Hopf algebras. The following cocycle $\omega'$ on $\widehat{M}$ has the advantage of being defined over $\Q$, so the fourth root of unity is dispensable:
$$\omega'(\varphi_1^{\alpha_1}\varphi_2^{\alpha_2},\varphi_1^{\beta_1}\varphi_2^{\beta_2}) = (-1)^{\alpha_1\beta_2}, \qquad \alpha_i,\beta_j \in \{0,1\}.$$
Put $L=K(\xi)$. Both cocycles are cohomologous over $L$: defining $q:\widehat{M} \rightarrow L^{\times}$ by
$$\varepsilon \mapsto 1, \quad \varphi_1 \mapsto 1, \quad \varphi_2 \mapsto 1, \quad \varphi_1\varphi_2 \mapsto \xi,$$
we have $\omega=\omega' \partial(q)$. Let $J'$ be the twist for $K\hspace{-1pt}A_5$ arising from $\omega'$. Then, $(LA_5)_J$ and $(LA_5)_{J'}$ are isomorphic. Now, if $X$ is a Hopf order of $(K\hspace{-1pt}A_5)_{J'}$ over $\Oint_K$, then $X \otimes_{\Oint_K} \hspace{-2pt}\Oint_L$ is a Hopf order of $(LA_5)_{J}$ over $\Oint_L$.
\end{remark}

We are now in a position to formulate and prove our first main result:

\begin{theorem}\label{main1}
Let $K$ be a number field containing a primitive $4$th root of unity. Let $R\subset K$ be a Dedekind domain such that $\Oint_K \subseteq R$. For $n \geq 5$ consider the twist $J$ for $K\hspace{-1pt}A_n$ given in \eqref{cocAn}. If $(K\hspace{-1pt}A_n)_J$ admits a Hopf order over $R$, then $\frac{1}{2}\in R$. \par \smallskip

In consequence, $(K\hspace{-1pt}A_n)_J$ does not admit a Hopf order over any number ring.
\end{theorem}

\pf Since $(K\hspace{-1pt}A_5)_J$ is a Hopf subalgebra of $(K\hspace{-1pt}A_n)_J$, it suffices to prove the statement for $(K\hspace{-1pt}A_5)_J$ in virtue of Proposition \ref{subsquo}(iii). Let $X$ be a Hopf order of $(K\hspace{-1pt}A_5)_J$ over $R$. The strategy of the proof consists in constructing $\nu \in X^{\star}$ and $x \in X$ such that $\nu(x)=\frac{27}{4}$; a number lying in $(\frac{1}{4}\Ow_K) \backslash \Ow_K$. By definition of $X^{\star}$, we have $\nu(x) \in R$. Thus $\frac{27}{4} \in R$ and from this it\hspace{0.6pt} follows\hspace{0.6pt} that\hspace{0.65pt} $\frac{1}{2} \in R$. The elements $\nu$ and $x$ will be constructed from characters and cocharacters of $(K\hspace{-1pt}A_5)_J$. \par \smallskip

Put $\tau=(12345)$. In Proposition \ref{A5Jdecomp}(ii) we saw that $\mu:=4e_{\varepsilon}\tau e_{\varepsilon}$ is a cocharacter of $(K\hspace{-1pt}A_5)_J$. By Proposition \ref{character}, we have $\mu \in X$. Then, $\Delta_J(\mu) \in X \otimes_R X$ because $X$ is a Hopf order.
We know from Proposition \ref{A5Jdecomp}(iii) that
\begin{equation}\label{comu}
\Delta_J(\mu)=\frac{1}{4}\sum_{\sigma \in M \tau M} \sigma \otimes \sigma.
\end{equation}
A direct calculation shows that
\begin{equation}\label{coset}
\hspace{-0.45cm}\begin{array}{l}
M\tau M \hspace{-0.5pt} = \hspace{-0.5pt} \big\{(12345),(14325),(15)(24),(135),(14532),(12534),(153),(24)(35),\vspace{3pt} \\
\hspace{1.746cm} (13)(45),(254),(15432),(12354),(245),(13)(25),(15234), (14352)\big\}.
\end{array}\hspace{-0.15cm}
\end{equation}
\par \smallskip

We next work with several characters of $(K\hspace{-1pt}A_5)_J$. The twist operation does not alter the algebra structure of the Hopf algebra, so the characters of $(K\hspace{-1pt}A_5)_J$ are those of $K\hspace{-1pt}A_5$. Recall from \cite[Section 18.6, p. 156]{Se} the character table of $A_5$:
$$\begin{array}{r|ccccc}
         &  id  &  (12)(34)  &  (123)  &  (12345)  &  (13524) \\
\hline
\chi_1   &   1  &     1     &    1    &      1    &     1     \\
\chi_2   &   3  &    \hspace{-8pt}-1     &    0    &  \frac{1+\sqrt{5}}{2} &  \frac{1-\sqrt{5}}{2} \\
\chi_3   &   3  &    \hspace{-8pt}-1     &    0    &  \frac{1-\sqrt{5}}{2} &  \frac{1+\sqrt{5}}{2} \\
\chi_4   &   4  &     0     &    1    &     \hspace{-8pt}-1    &     \hspace{-8pt}-1     \\
\chi_5   &   5  &     1     &   \hspace{-8pt}-1    &      0    &      0     \\
\end{array}$$
The conjugacy class of $(12)(34)$ contains all double transpositions of $A_5$. All $3$-cycles belong to the conjugacy class of $(123)$. \par \smallskip

Consider the character $$\chi=\chi_1+\chi_4+2\chi_5.$$
Proposition \ref{character} yields $\chi \in X^{\star}$. Observe that $\chi$ gives $3$ when evaluated at any double transposition and $\chi$ vanishes at the conjugacy classes of either $(123), (12345)$ or $(13524)$. By construction, the following element $y$ belongs to $X$:
$$\begin{array}{rl}
y & \hspace{-1.5mm}:= \hspace{1.5pt}(\chi \otimes_R Id_X)\Delta_J(\mu) \vspace{5pt} \\
  & \hspace{-1.5mm}\hspace{-3pt}\overset{\eqref{comu}}{=} \hspace{1.5pt}{\displaystyle \frac{1}{4}\sum_{\sigma \in M \tau M} \chi(\sigma) \sigma} \vspace{5pt} \\
  & \hspace{-1.5mm}\hspace{-3pt}\overset{\eqref{coset}}{=} \hspace{1.5pt}{\displaystyle \frac{3}{4}\big((15)(24)+(24)(35)+(13)(45)+(13)(25)\big).}
\end{array}$$
Since $X$ is a Hopf order, $y^2$ belongs to $X$ as well. A computation gives:
$$\begin{array}{l}
{\displaystyle y^2 = \frac{9}{16}\big(4 \cdot id + (15243)+(13425)+(15423)+(13245)+(135)+(153)} \vspace{1.6pt} \\
{\displaystyle \phantom{a} \hspace{1.3cm} + (452)+(425)+(14253)+(13524)+(12453)+(13542)\big).}
\end{array}$$

Finally, applying $\chi_4$ to $y^2$ we get the desired value,
$$\chi_4(y^2)=\frac{9}{16}\big(4 \cdot 4 + 4 \cdot 1 + 8 \cdot (-1)\big) = \frac{27}{4};$$
which resides in $R$ as $\chi_4 \in X^{\star}$.
\epf

\begin{remark}
In view of Proposition \ref{subsquo}(iii), the thesis of Theorem \ref{main1} holds for $(K\hspace{-0.9pt}G)_J$, where $G$ is any finite group containing $A_5$. In particular, it holds for the symmetric group $S_n$ and for the projective special linear group $PSL(2,5^r)$, with $r \geq 1$.
\end{remark}

The same argument exhibited in the proof of \cite[Corollary 2.4]{CM} establishes:

\begin{corollary}\label{complexif}
For $n \geq 5$ the complex semisimple Hopf algebra $(\Co A_n)_J$ does not admit a Hopf order over any number ring.
\end{corollary}

\begin{remark}
The proof of the non-existence of Hopf orders for $(\Co A_n)_J$ is simpler than for $\mathcal{B}_{p,q}(\zeta)$, the first family for which we proved this property in \cite[Theorem 2.3]{CM}. Here we use neither the character supported algebra, \cite[Definition 1.5 and Proposition 1.6]{CM}, nor the preservation of Hopf orders under surjective Hopf algebra maps. Compare also the computations here with those in the intricate proof of \cite[Proposition 2.1]{CM}, which is the support of the main result.
\end{remark}

\section{Twisting of the symmetric group}

In this section we study the existence of Hopf orders for a twist of the symmetric group that was introduced by Bichon in \cite{B} and further discussed by Galindo and Natale in \cite[Section 3]{GN}. \par \smallskip

Take $n \in \Na$ even such that $n \geq 2$. For $1 \leq i \leq n$ let $t_i$ denote the transposition $(2i-1\ 2i)$ in $S_{2n}$. Consider the abelian subgroup $M$ of $S_{2n}$ generated by the elements\vspace{-1pt} $t_i$ for $1 \leq i \leq n$. It is isomorphic to $\Z_2^n.$ The character group $\widehat{M}$ is generated by $\varphi_i$,  for $1 \leq i \leq n,$ given by $\varphi_i(t_j)=(-1)^{\delta_{ij}}$.  \par \smallskip

The bicharacter $\omega:\widehat{M} \times \widehat{M} \rightarrow K^{\times}$ defined as
$$\omega(\varphi_i,\varphi_j)=\left\{\hspace{-4pt}\begin{array}{rl}
-1 & \textrm{if}\hspace{2pt} i<j, \vspace{3pt} \\
 1 & \textrm{if}\hspace{2pt} i \geq j,
\end{array}\right.$$
is a $2$-cocycle on $\widehat{M}$. The general formula for $\omega$ is
$$\omega(\varphi_1^{\alpha_1}\ldots \varphi_n^{\alpha_n},\varphi_1^{\beta_1}\ldots \varphi_n^{\beta_n}) = (-1)^{\sum_{i<j}\alpha_i\beta_j}, \qquad \alpha_i,\beta_j \in \{0,1\}.$$
For $(\alpha_1,\ldots,\alpha_n) \in \{0,1\}^n$, consider the idempotent of $K\hspace{-1pt}M$ given by
$$e_{(\alpha_1,\ldots,\alpha_n)} = \frac{1}{\vert M \vert} \sum_{m \in M} \langle \varphi_1^{\alpha_1}\ldots \varphi_n^{\alpha_n}, m \rangle m.$$

Take the twist
\begin{equation}\label{cocSn}
J=\sum_{\phi,\psi \in \widehat{M}} \omega(\phi,\psi) e_{\phi} \otimes e_{\psi}
\end{equation}
for $K\hspace{-1pt}M$ afforded by $\omega$.
\par \smallskip

The analysis of the existence of Hopf orders for $(K\hspace{-1pt}S_{2n})_J$ needs to distinguish two cases: $n=2$ and $n>2$. For $n=2$ a Hopf order exists. We address this case first. We will use several arguments that will appear in the general case too, although they take a simpler form here.

\subsection{The case $n=2$.}
The subgroup $M$ of $S_4$ is now generated\vspace{-1pt} by $t_1=(12)$ and $t_2=(34)$. The cocycle $\omega$ on $\widehat{M}$ is given by the following table:
\begin{equation}\label{omegaS4}
\begin{array}{c|rrrr}
                     &  \varepsilon  & \varphi_1 & \varphi_2 & \varphi_1\varphi_2 \\
\hline
\varepsilon          &  1  &  1  &   1  &   1   \\
\varphi_1        	 &  1  &  1  &  -1  &  -1   \\
\varphi_2            &  1  &  1  &   1  &   1  \\
\varphi_1\varphi_2   &  1  &  1  &  -1  &  -1
\end{array}
\end{equation}

The orthogonal idempotents providing the Wedderburn decomposition of $K\hspace{-0.9pt}M$ are:
\begin{equation}\label{idemS4}
\begin{array}{lp{1cm}l}
{\displaystyle e_{(0,0)} = \frac{1}{4}(id+t_1+t_2+t_1t_2),}  &  & {\displaystyle e_{(1,0)} = \frac{1}{4}(id-t_1+t_2-t_1t_2),} \vspace{5pt} \\
{\displaystyle e_{(0,1)} = \frac{1}{4}(id+t_1-t_2-t_1t_2),}  &  & {\displaystyle e_{(1,1)} = \frac{1}{4}(id-t_1-t_2+t_1t_2).} \\
\end{array}
\end{equation}
\smallskip

One can easily check that $J$ reads as:
\begin{equation}\label{JS4}
J= \frac{1}{2}\big(id+(12)\big) \otimes id +  \frac{1}{2}\big(id-(12)\big) \otimes (34).
\end{equation}

Our aim in this part is to prove the following result:

\begin{proposition}\label{caseS4}
Let $K$ be a number field containing a primitive 4th root of unity $\xi$. The Hopf algebra $(K\hspace{-1pt}S_4)_J$ admits a Hopf order over $\Ow_K$. Moreover, it is unique.
\end{proposition}

\pf We first deal with the existence. Let $F$ denote the subgroup of $S_4$ generated by $(12)(34)$ and $(13)(24)$. Consider the following idempotents in $K\hspace{-1pt}F$:
$$x=\frac{1}{2}\big(id-(12)(34)\big) \quad \textrm{and}\quad y=\frac{1}{2}\big(id-(13)(24)\big).$$
Observe that $x \in K\hspace{-1pt}F \cap K\hspace{-1pt}M$ and $y \notin K\hspace{-1pt}M$. Let $X$ be the $\Ow_K$-submodule of $K\hspace{-1pt}S_4$ generated by the set
\begin{equation}\label{spanningset}
\{xy\sigma,\, (id-x)y\sigma,\, x(id-y)\sigma,\, (id-x)(id-y)\sigma \, :\, \sigma \in S_3\}.
\end{equation}
We check that $X$ is a Hopf order of $K\hspace{-1pt}S_4$ over $\Ow_K$. Notice that $X$ is an order of $K\hspace{-1pt}S_4$ as the generating set is a basis. It contains $id$ and is closed under the product. For the latter, bear in mind that $x$ and $y$ are idempotents and that conjugation by $\sigma$ permutes the elements $xy, (id-x)y, x(id-y),$ and $(id-x)(id-y)$. To see that $X$ is closed under the coproduct, use that $\sigma$ is a group-like element and the formulas:
$$\Delta(x)= id \otimes x + x \otimes (12)(34) \quad \textrm{and} \quad \Delta(y)= id \otimes y + y \otimes (13)(24).$$
That $\varepsilon(X) \subseteq \Ow_K$ and $\mathcal{S}(X)\subseteq X$ are clear.\par \smallskip

Next, we will obtain a Hopf order of $(K\hspace{-1pt}S_4)_J$ from $X$. If $J$ would belong to $X \otimes_{\Ow_K} X$, we could construct the twisted Hopf order $X_J$ in virtue of Proposition \ref{twistorders} (notice that $J=J^{-1}$). However, this is not true, as Remark \ref{JnotinXX} shows. Nevertheless, that condition can be achieved for a cohomologous twist. We replace the cocycle $\omega$ by a cohomologous one, $\kappa$, which is given by the following table:
$$\begin{array}{c|rrrr}
             &  \varepsilon & \varphi_1 & \varphi_2 & \varphi_1\varphi_2 \\
\hline
\varepsilon  &  1  &  1  &  1  &  1   \\
\varphi_1    &  1  &  -1  &  1  & -1  \\
\varphi_2    &  1  & -1  &  1  &  -1  \\
\varphi_1\varphi_2  &  1  & 1  &  1  &  1
\end{array}$$
Define $q:\widehat{M} \rightarrow K^{\times}$ by
\begin{equation}\label{cohomq}
\varepsilon \mapsto 1, \quad \varphi_1 \mapsto \xi, \quad \varphi_2 \mapsto -1, \quad \varphi_1\varphi_2 \mapsto \xi.
\end{equation}
Then, $\omega=\kappa \partial(q)$. Let $T$ be the twist for $K\hspace{-1pt}M$ afforded by $\kappa$. There is $v \in K\hspace{-1pt}M$ invertible such that the map
$$f:(K\hspace{-1pt}S_4)_T \rightarrow (K\hspace{-1pt}S_4)_J,\, h \mapsto vhv^{-1},$$
is a Hopf algebra isomorphism. An easy computation reveals that:
\begin{equation}\label{twistT}
T = (id-x)\otimes id + x \ot (12).
\end{equation}
Hence, $T \in X\ot_{\Ow_K} \hspace{-2pt}X$. Since $T=T^{-1}$, by Proposition \ref{twistorders}, $X_T$ is a Hopf order of $(K\hspace{-1pt}S_4)_T$, and so $f(X_T)$ is a Hopf order of $(K\hspace{-1pt}S_4)_J$. \par \medskip

We now tackle the uniqueness. We work in $(K\hspace{-1pt}S_4)_T$. Let $Y$ be a Hopf order of $(K\hspace{-1pt}S_4)_T$ over $\Ow_K$. We will first show that $X \subseteq Y$. By manipulating characters and cocharacters of $(K\hspace{-1pt}S_4)_J$ we will derive that the generating set in \eqref{spanningset} is included in $Y$.  \par \smallskip

There are $3$ double cosets of $M$ in $S_4$, with cardinalities $4,4$ and $16$; namely:
\begin{flushleft}
$\begin{array}{rl}
         M & \hspace{-2mm}= \hspace{1.5mm}\{id,(12),(34),(12)(34)\}, \vspace{5pt} \\
M(13)(24)M & \hspace{-2mm}= \hspace{1.5mm}\{(13)(24),(14)(23),(1324),(1423)\}  \vspace{3pt} \\
           & \hspace{-2mm}= \hspace{1.5mm}M(13)(24), \vspace{5pt} \\
M(123)M    & \hspace{-2mm}=  \hspace{1.5mm}\left\{(123),(13),(1234),(134),(23),(132),(234),(1342),\right. \vspace{3pt} \\
           & \phantom{=}\hspace{6pt} \left. (143),(14),(1432),(124),(24),(243),(142),(1243)\right\}.
\end{array}$
\end{flushleft}

Set $\tau=(123)$. Observe that $\omega$ is non-degenerate and $M \cap (\tau M \tau^{-1})=\{id\}$. By Proposition \ref{simple}, $K(M\tau M)$ is isomorphic to a matrix coalgebra of size $4$ and the irreducible cocharacter afforded by it is $\mu:=4e_{\varepsilon}\tau e_{\varepsilon}$. Write $\overline{Y}$ for the Hopf order $f(Y)$ of $(K\hspace{-1pt}S_4)_J$. By Proposition \ref{character}, we have $\mu \in \overline{Y}$. Then, $\Delta_J(\mu) \in \overline{Y} \otimes_{\Ow_K} \overline{Y}$. We compute:
\begin{align}
\Delta_J(\mu) & \hspace{2pt}\overset{\eqref{eqcopr}}{=} {\displaystyle \sum_{i,j,k,l} \frac{\omega(\varphi_1^i\varphi_2^j,\varphi_1^i\varphi_2^j)}{\omega(\varphi_1^k\varphi_2^l,\varphi_1^k\varphi_2^l)} e_{(i,j)}\tau e_{(k,l)}\ot e_{(i,j)}\tau e_{(k,l)}} \vspace{7pt} \notag \\
              & \hspace{2pt}\overset{\eqref{omegaS4}}{=} {\displaystyle 4\sum_{i,j,k,l}(-1)^{ij+kl}e_{(i,j)}\tau e_{(k,l)}\ot e_{(i,j)}\tau e_{(k,l)}.} \label{expforr}
\end{align}

Recall from \cite[Section 5.8, p. 43]{Se} the character table of $S_4$:
\begin{equation}\label{ctabS4}
\begin{array}{r|ccccc}
         &  id  &  (12) &  (12)(34)  &  (123)  &  (1234) \\
\hline
\chi_1   &   1  &     1     &    1   &   1   &  1     \\
\chi_2   &   1  &  \hspace{-8pt}-1   &   1   &  1 &  \hspace{-8pt}-1   \\
\chi_3   &   2  &    0      &    2   &  \hspace{-8pt}-1   &  0 \\
\chi_4   &   3  &     1     &   \hspace{-8pt}-1  &   0   &  \hspace{-8pt}-1 \\
\chi_5   &   3  &  \hspace{-8pt}-1   &   \hspace{-8pt}-1    &      0    &  1 \\
\end{array}
\end{equation}
We now calculate the element
$$r:=(\chi_4 \ot Id_{K\hspace{-1pt}S_4})\Delta_J(\mu),$$
which will belong to $\overline{Y}$. First, we use that $\chi_4$ is a character in the following way:
$$\chi_4(e_{(\alpha_1,\alpha_2)}\tau e_{(\beta_1,\beta_2)}) = \chi_4(\tau e_{(\beta_1,\beta_2)}e_{(\alpha_1,\alpha_2)})= \chi_4(\tau e_{(\alpha_1,\alpha_2)})\delta_{\alpha_1,\beta_1}\delta_{\alpha_2,\beta_2}.$$
Secondly, we observe that:
$$\chi_4(\tau) = 0, \hspace{8pt} \chi_4(\tau t_1)=1, \hspace{8pt} \chi_4(\tau t_2)=-1,\hspace{8pt} \textrm{and} \hspace{8pt} \chi_4(\tau t_1t_2)=0.$$
Using the Formula \ref{idemS4}, which express $e_{(i,j)}$ in terms of $t_1$ and $t_2$, we get:
$$\chi_4(\tau e_{(0,0)})=0, \hspace{8pt} \chi_4(\tau e_{(0,1)})=\frac{1}{2}, \hspace{8pt} \chi_4(\tau e_{(1,0)})=-\frac{1}{2}, \hspace{8pt} \textrm{and} \hspace{7pt} \chi_4(\tau e_{(1,1)})=0.$$
We compute $r$ bearing this in mind. We have:
\begin{equation}\label{elementr}
\begin{array}{rl}
r & \hspace{-3pt}\overset{\eqref{expforr}}{=} {\displaystyle 4 \sum_{i,j,k,l}(-1)^{ij+kl} \chi_4 \big(e_{(i,j)}\tau e_{(k,l)}\big) e_{(i,j)}\tau e_{(k,l)}} \vspace{5pt} \\
  & \hspace{3pt}= \hspace{2pt}2\big(e_{(0,1)}\tau e_{(0,1)}-e_{(1,0)}\tau e_{(1,0)}\big).
\end{array}
\end{equation}
We next write $r$ as a linear combination of elements in $S_4$. A calculation yields:
$$\begin{array}{rl}
e_{(0,1)}\tau e_{(0,1)} & \hspace{-2mm} = {\displaystyle \frac{1}{16}\big((123)+(13)-(1234)-(134)+(23)+(132)-(234)-(1342)} \\
 & {\displaystyle \hspace{7mm}-(1243)-(143)+(124)+(14)-(243)-(1432)+(24)+(142)\big).} \vspace{5pt} \\
e_{(1,0)}\tau e_{(1,0)} & \hspace{-2mm} = {\displaystyle \frac{1}{16}\big((123)-(13)+(1234)-(134)-(23)+(132)-(234)+(1342)} \\
 & {\displaystyle \hspace{7mm}+(1243)-(143)+(124)-(14)-(243)+(1432)-(24)+(142)\big).} \vspace{5pt} \\
\end{array}$$
We finally obtain:
$$r=\frac{1}{4}\big((13)-(1234)+(23)-(1342)-(1243)+(14)-(1432)+(24)\big).$$

With $r$ in our hands, we compute $r^2$ by a direct and tedious calculation (or with the help of a computer). We get:
$$r^2 = \frac{1}{2}\big(id-(12)(34)\big)=x.$$
Thus, $x$ and $id-x$ belong to $\overline{Y}$. Observe that $f\vert_ M = Id_M$. So $x$ and $id-x$ belong to $Y$. On the other hand, $(12)$ is a group-like element of $(K\hspace{-1pt}S_4)_T$. Proposition \ref{character} entails that $(12) \in Y$. The Formula \ref{twistT} for $T$ implies that $T \in Y \otimes_{\Ow_K} \hspace{-2pt} Y$. Since $T=T^{-1}$, the twisted order $Y_{T^{-1}}$ is a Hopf order of $K\hspace{-1pt}S_4$ by Proposition \ref{twistorders}. Using Proposition \ref{character} once more, we have $S_3 \subset Y$. Consequently, $y=(23)x(23)$ belongs to $Y$. This proves that $X \subseteq Y$. \par \smallskip

Recall that $F$ is the (normal) subgroup of $S_4$ generated by $(12)(34)$ and $(13)(24)$. For simplicity, we denote the idempotents of $K\hspace{-1pt}F$ by $\epsilon_{ij}$, for $i,j=0,1$; where $(12)(34)\epsilon_{ij} = (-1)^i \epsilon_{ij}$ and $(13)(24) \epsilon_{ij} = (-1)^j \epsilon_{ij}$. Notice that
$$\epsilon_{00}=(id-x)(id-y), \quad \epsilon_{10}=x(id-y), \quad \epsilon_{01}=(id-x)y, \quad \textrm{and} \quad \epsilon_{11}=xy.$$
Then, $\epsilon_{ij}$ lies in $Y$ for all $i,j$ as $X \subseteq Y$. \par \smallskip

We next check that $Y \subseteq X$. Pick $w \in Y$. Since $S_4=FS_3$, we can write $w$ in the form $w=\sum_{i,j} \epsilon_{ij} a_{ij}$ with $a_{ij}\in K\hspace{-1pt}S_3$. In turn, write $a_{ij}=\sum_{\sigma \in S_3} \lambda_{ij\sigma} \sigma$ with $\lambda_{ij\sigma} \in K$. We will prove that $\lambda_{ij\sigma} \in \Ow_K$ and this will yield $Y \subseteq X$. Consider the element $\Delta^2(w)$ (coproduct in $K\hspace{-1pt}S_4$). It is:
$$\Delta^2(w)= \sum_{i,j} \sum_{i_1,i_2,j_1,j_2} (\epsilon_{i_1,j_1} \otimes \epsilon_{i_2,j_2}\otimes \epsilon_{i_1+i_2+i,\hspace{1pt} j_1+j_2+j})\Delta^2(a_{ij}).$$
Multiply it by $\epsilon_{01}\otimes \epsilon_{10} \otimes \epsilon_{11}$ from both sides. We get:
$$\begin{array}{ll}
(\epsilon_{01}\otimes \epsilon_{10} \otimes \epsilon_{11})\Delta^2(w)(\epsilon_{01}\otimes \epsilon_{10} \otimes \epsilon_{11}) & = (\epsilon_{01}\otimes \epsilon_{10}\otimes \epsilon_{11}) \Delta^2(a_{00}) (\epsilon_{01}\otimes \epsilon_{10}\otimes \epsilon_{11}) \vspace{5pt} \\
  & {\displaystyle = \sum_{\sigma \in S_3} \lambda_{00\sigma} \epsilon_{01}\sigma \epsilon_{01} \otimes \epsilon_{10}\sigma \epsilon_{10} \otimes \epsilon_{11}\sigma \epsilon_{11}.}
\end{array}$$
Observe that the action of $S_3$ by conjugation on $\{\epsilon_{01},\epsilon_{10},\epsilon_{11}\}$ satisfies the following property: if $\sigma \in S_3$ fixes any two distinct elements, then $\sigma=id$. Taking into account this and that these idempotents are orthogonal, we have:
$$(\epsilon_{01}\otimes \epsilon_{10} \otimes \epsilon_{11})\Delta^2(w)(\epsilon_{01}\otimes \epsilon_{10} \otimes \epsilon_{11}) = \lambda_{00id} \hspace{1pt} \epsilon_{01}\otimes \epsilon_{10} \otimes \epsilon_{11}.$$
The left-hand side term belongs to $Y \otimes_{\Ow_K} \hspace{-2pt}Y \otimes_{\Ow_K} \hspace{-2pt}Y$. So, the right-hand side term satisfies a monic polynomial with coefficients in $\Ow_K$. Since $\epsilon_{01} \otimes \epsilon_{10} \otimes \epsilon_{11}$ is idempotent, it must hold that $\lambda_{00id} \in \Ow_K$. Arguing as before with $\sigma^{-1} w$ we obtain that $\lambda_{00\sigma} \in \Ow_K$ for all $\sigma \in S_3$. \par \smallskip

Finally, by multiplying $\Delta^2(w)$ from both sides with other triples of idempotents we get in a similar fashion that $\lambda_{ij\sigma} \in \Ow_K$ for all $i,j,$ and $\sigma.$ Concretely, use $\epsilon_{10} \otimes \epsilon_{11} \otimes \epsilon_{11}$ for $a_{10}$; $\epsilon_{01} \otimes \epsilon_{11} \otimes \epsilon_{11}$ for $a_{01};$ and $\epsilon_{01} \otimes \epsilon_{11} \otimes \epsilon_{01}$ for $a_{11}.$
\epf \vspace{-2mm}

\begin{remark}\label{JnotinXX}
The following argument shows that $J \notin X \otimes_{\Ow_K} \hspace{-2pt}X$. Assume the contrary. Since $X$ is a Hopf order of $K\hspace{-1pt}S_4$ over $\Ow_K$, Proposition \ref{character} gives that $\chi_3 \in X^{\star}$ and $S_4 \subset X$. Then, $d:=J\big((123) \otimes (134)\big)$ belongs to $X \otimes_{\Ow_K} \hspace{-2pt}X$ and, consequently, $(\chi_3 \otimes \chi_3)(d)$ lies in $\Ow_K$. We compute:
$$\begin{array}{rl}
(\chi_3 \otimes \chi_3)(d) & \overset{\eqref{JS4}}{=} {\displaystyle \frac{1}{2}\chi_3\big((123)+(23)\big)\chi_3\big((134)\big) +  \frac{1}{2}\chi_3\big((123)-(23)\big)\chi_3\big((14)\big)} \vspace{6pt} \\
  & \hspace{-2pt}\overset{\eqref{ctabS4}}{=} {\displaystyle \frac{1}{2}.}
\end{array}$$
And this is a contradiction.
\end{remark}

\begin{remark}
The root of unity $\zeta$ is not needed to construct the Hopf order $X_T$ of $(K\hspace{-1pt}S_4)_T$, but only to claim that $\omega$ and $\kappa$ are cohomologous (definition of $q$ in \eqref{cohomq}); that is, to claim that $(K\hspace{-1pt}S_4)_J$ and $(K\hspace{-1pt}S_4)_T$ are isomorphic. One could start with $\kappa$ instead of $\omega$ and all statements would hold without requiring the existence of $\zeta$. We used $\omega$ just to keep the definition of \cite{GN}. \par \smallskip

The Hopf order of $(K\hspace{-1pt}S_4)_J$ can explicitly be described by carrying $X$ through the isomorphism $f:(K\hspace{-1pt}S_4)_T \rightarrow (K\hspace{-1pt}S_4)_J, \, h \mapsto vhv^{-1}.$ The only disadvantage is the cumbersome form of the image of $(123)$. In view of \eqref{cohomq}, we have:
$$v = e_{(0,0)}+\xi e_{(1,0)}-e_{(0,1)}+\xi e_{(1,1)} =\frac{1}{2}\big(\xi id-\xi (12)+(34)+(12)(34)\big).$$
The unique Hopf order of $(K\hspace{-1pt}S_4)_J$ over $\Ow_K$ is the $\Ow_K$-subalgebra of $(K\hspace{-1pt}S_4)_J$ generated by the elements $x,y,(12),$ and
$$\begin{array}{l}
{\displaystyle \frac{1}{4}\big((123)-(13)+\xi(1234)+\xi(134)-(23)+(132)-\xi(234)-\xi(1342)} \vspace{3pt} \\
\hspace{1cm}  {\displaystyle -\xi(1243)+\xi(143)+(124)+(14)-\xi(243)+\xi(1432)+(24)+(142)\big).}
\end{array}$$
\end{remark}
\vspace{1mm}

\begin{remark}
The following interpretation of the Hopf order $X$ of $K\hspace{-1pt}S_4$ and its consequence was explained to us by the referee. Recall that $F$ is the subgroup of $S_4$ generated by $(12)(34)$ and $(13)(24)$.\vspace{-1pt} The idempotents $x$ and $y$ generate the unique maximal order $A$ of $K\hspace{-1pt} F$. This is the dual of the minimal Hopf order $\Oint_K \widehat{F}$ of $K\hspace{-1pt} \widehat{F}$. As $F$ is the unique nontrivial normal subgroup of $S_4$ that is a $p$-group, by \cite[Proposition 3.1 and Corollary 3.8]{L2} and \cite[Corollary 17.4]{C}, $X=A(\Oint_K S_4)$ must be the maximal Hopf order of $K\hspace{-1pt} S_4$. There are other (smaller) Hopf orders $B$ in $K\hspace{-1pt} S_4$, for example $\Oint_K S_4$. The uniqueness of $X_T$ in $(K\hspace{-1pt}S_4)_T$ shows that none of them can satisfy $T \in B \otimes_{\Oint_K} \hspace{-2pt}B$.
\end{remark}
\vspace{1mm}

This concludes our discussion of the case $n=2$. We next address the general case.

\subsection{The case $n \geq 4$.} Our second main result is formulated as follows:

\begin{theorem}\label{main2}
Let $K$ be a number field and $R\subset K$ be a Dedekind domain such that $\Oint_K \subseteq R$. For $n \geq 4$ even, consider the twist $J$ for $K\hspace{-1pt}S_{2n}$ given in \eqref{cocSn}. If $(K\hspace{-1pt}S_{2n})_J$ admits a Hopf order over $R$, then $\frac{1}{2}\in R$. \par \smallskip

As a consequence, $(K\hspace{-1pt}S_{2n})_J$ does not admit a Hopf order over any number ring.
\end{theorem}

\pf
Assume that $X$ is a Hopf order of $(K\hspace{-1pt}S_{2n})_J$ over $R$. The proof is organized in several steps and proceeds by various reductions: \par \medskip

{\it Step 1. Reduction to the case $n=4$.} Decompose $M$ as $M=PQ$, with $P=\langle (12),\ldots, (78) \rangle$ and $Q=\langle (9\hspace{2pt} 10),\ldots, (2n-1\hspace{2.5pt} 2n)\rangle$. Consider the subgroup $S_8Q$ of $S_{2n}$. Since $J$ is supported on $M$ and $M \subset S_8Q$, we can construct the twisted Hopf algebra $K(S_8Q)_J$, which is a Hopf subalgebra of $(K\hspace{-1pt}S_{2n})_J$. By Proposition \ref{subsquo}(iii), $X \cap (K(S_8Q)_J)$ is a Hopf order of $K(S_8Q)_J$ over $R$. \par \smallskip

The subgroup $S_8$ commutes with $Q$ and $S_8 \cap Q=\{id\}$. Identify $S_8Q$ with the direct product group $S_8 \times Q$. Projecting on the first factor induces a surjective\vspace{-0.5pt} Hopf algebra map $\pi:K(S_8Q)_J \to (KS_8)_{\tilde{J}}$, where $\tilde{J}=(\pi \otimes \pi)(J)$. Identify\vspace{-1pt} now $\widehat{M}$ with $\widehat{P} \times \widehat{Q}$. For $\phi \in \widehat{M}$, the idempotent $e_{\phi} \in K\hspace{-1pt}M$ equals $e_{\phi_1}e_{\phi_2}$ with $\phi_1 \in \widehat{P}$ and $\phi_2 \in \widehat{Q}$ such that $\phi=(\phi_1,\phi_2)$. Observe that $\pi(e_{\phi})=e_{\phi_1}\delta_{\varepsilon,\phi_2}$. Then, $\tilde{J}$ is precisely\vspace{0.5pt} the twist $J$ of $S_8$ afforded by $\omega \vert_{\widehat{P} \times \widehat{P}}$. In view of Proposition\vspace{0.5pt} \ref{subsquo}(iv), $\pi(X \cap (K(S_8Q)_J))$ is a Hopf order of $(K\hspace{-1pt}S_8)_J$ over $R$. \par \smallskip

We restrict our attention to $(K\hspace{-1pt}S_8)_J$. Set $H=(K\hspace{-1pt}S_8)_J$ and suppose that $X$ is now a Hopf order of $H$ over $R$. \par \medskip

{\it Step 2. Finding several elements in $X$.} In this step we will construct several elements in $X$ by manipulation of characters and cocharacters of $H$ and applying repeatedly Proposition \ref{character}. \par \smallskip

Recall that $M$ is generated by $t_i=(2i-1\hspace{2.5pt} 2i)$ with $i=1,2,3,4.$ Decompose again $M$ as $M=PQ$, where this time $P=\langle (12),(34) \rangle$ and $Q=\langle (56),(78)\rangle$. Notice that $S_4$ commutes with $Q$ and $S_4 \cap Q=\{id\}$.
Consider the subgroup $G=S_4Q$. Since $M \subset G$, we can construct the twisted Hopf algebra $(K\hspace{-0.9pt}G)_J$, which is a Hopf subalgebra of $H$. Put $A=(K\hspace{-0.5pt}G)_J$. By Proposition \ref{subsquo}(iii), $X \cap A$ is a Hopf order of $A$ over $R$. We now focus on $A$ and the Hopf order $X \cap A$. \par \smallskip

Take $\tau=(123)$. We will find a cocharacter of $A$ arising from the subcoalgebra $K(M\tau M)$. Observe that $M \cap (\tau M \tau^{-1}) = \langle t_3,t_4\rangle$. Consider the following subgroup of $\widehat{M} \times \widehat{M}$:
$$\hspace{1cm} N=\big\{(\phi,\psi) \in \widehat{M}\times\widehat{M} \ : \ \psi(m)=\phi(\tau m \tau^{-1}) \hspace{7pt} \forall m \in M\cap (\tau M \tau^{-1}) \big\}.$$
We know that
$$\vert N \vert = \vert M \tau M \vert = \frac{|M|^2}{\vert M \cap (\tau M \tau^{-1}) \vert }= \frac{256}{4}=64.$$
A direct computation shows that $N$ is generated by the following six elements:
$$(\varepsilon,\varphi_1),(\varepsilon,\varphi_2),(\varphi_1,\varphi_1), (\varphi_2,\varphi_2), (\varphi_3,\varphi_3), \ \textrm{and}\ (\varphi_4,\varphi_4).$$
We saw in the proof of Proposition \ref{decomp} that $\{e_{\phi}\tau e_{\psi}\}_{(\phi,\psi) \in N}$ is a basis of $K(M\tau M)$ and that $K(M\tau M)^* \simeq K^{(\omega,\omega^{-1})\vert_N}[N]$. Recall also from there that $\{u_{(\phi,\psi)}\}_{(\phi,\psi)\in N}$ denotes the dual basis of the previous one. Since $N$ is abelian, the center of $K^{(\omega,\omega^{-1})\vert_N}[N]$ is spanned by the elements $u_{(\phi,\psi)}$ such that
$$\frac{\omega(\phi,\phi')}{\omega(\psi,\psi')}=\frac{\omega(\phi',\phi)}{\omega(\psi',\psi)}, \qquad \forall (\phi',\psi') \in N.$$
The set of pairs obeying this condition is: $$Z:=\{(\varepsilon,\varepsilon),(\varphi_1\varphi_2\varphi_3,\varphi_1\varphi_2\varphi_3),(\varphi_3\varphi_4,\varphi_3\varphi_4),(\varphi_1\varphi_2\varphi_4,\varphi_1\varphi_2\varphi_4)
\}.$$
Thus, the center is generated, as an algebra, by $u_{(\varphi_1\varphi_2\varphi_3,\varphi_1\varphi_2\varphi_3)}$ and $u_{(\varphi_3\varphi_4,\varphi_3\varphi_4)}$. It has dimension $4$. All irreducible representations of $K^{(\omega,\omega^{-1})\vert_N}[N]$ are of dimension $4$. Consider the primitive central idempotent
$$c=\frac{1}{4}\big(u_{(\varepsilon,\varepsilon)}+u_{(\varphi_1\varphi_2\varphi_3,\varphi_1\varphi_2\varphi_3)}+u_{(\varphi_3\varphi_4,\varphi_3\varphi_4)}+u_{(\varphi_1\varphi_2\varphi_4,\varphi_1\varphi_2\varphi_4)}\big)
$$
and the representation $V=K^{(\omega,\omega^{-1})\vert_N}[N]c$. We know that $V$ is the direct sum of four copies of an irreducible representation. We next compute the corresponding irreducible character. Choose $\Gamma \subset N$ such that $\{u_{(\phi,\psi)}c:(\phi,\psi) \in \Gamma\}$ is a basis of $V$. Take $\gamma \in \Gamma$ and $\nu \in N$. If $\nu \in Z$, then $u_{\nu}u_{\gamma}c=u_{\gamma}u_{\nu}c=u_{\gamma}c$. If $\nu \notin Z$, then $u_{\nu}u_{\gamma}c \neq \beta u_{\gamma}c$ for every $\beta \in K$. Hence, the character of $K^{(\omega,\omega^{-1})\vert_N}[N]$ associated to $V$ maps $u_{\nu}$ to $16$ if $\nu \in Z$ and to $0$ otherwise. The irreducible character $\Theta$ of $K^{(\omega,\omega^{-1})\vert_N}[N]$ afforded by $V$ must be:
$$\Theta(u_{\nu})= \left\{\hspace{-4pt}\begin{array}{rl}
4 & \textrm{if}\hspace{6pt} \nu \in Z, \vspace{3pt} \\
0 & \textrm{if}\hspace{6pt} \nu \notin Z.
\end{array}\right.$$

By dualizing, we get the following irreducible cocharacter of $K(M\tau M)$:
$$\Psi = 4\big(e_{(0,0,0,0)}\tau e_{(0,0,0,0)} + e_{(1,1,1,0)}\tau e_{(1,1,1,0)} + e_{(0,0,1,1)}\tau e_{(0,0,1,1)} + e_{(1,1,0,1)}\tau e_{(1,1,0,1)}\big).$$

Identify $A$ with $K\hspace{-1pt}S_4 \otimes KQ$ as algebras. Consider the character $\theta=\chi_4 \otimes \varepsilon_Q$ on $A$, where $\chi_4$ is the character of $S_4$ given in Table \ref{ctabS4} and $\varepsilon_Q$ is the trivial character on $Q$. We next calculate the element
$$x:=(Id_{A} \ot \theta)\Delta_J(\Psi).$$
We first compute $\Delta_J(\Psi)$:
$$\begin{array}{rl}
\Delta_J(\Psi) \overset{\text{\eqref{eqcopr}}}{=} & \hspace{-2mm}{\displaystyle 4\scalebox{1.2}{$\Bigg($}\sum_{\substack{\alpha_1,\alpha_2,\alpha_3,\alpha_4 \\ \beta_1,\beta_2,\beta_3,\beta_4}}
\frac{\omega((\alpha_1,\alpha_2,\alpha_3,\alpha_4),(\alpha_1,\alpha_2,\alpha_3,\alpha_4))}{\omega((\beta_1,\beta_2,\beta_3,\beta_4),
(\beta_1,\beta_2,\beta_3,\beta_4))}} \vspace{3pt} \\
             &  \hfill e_{(\alpha_1,\alpha_2,\alpha_3,\alpha_4)} \tau e_{(\beta_1,\beta_2,\beta_3,\beta_4)}\ot   e_{(\alpha_1,\alpha_2,\alpha_3,\alpha_4)} \tau e_{(\beta_1,\beta_2,\beta_3,\beta_4)}
\end{array}\vspace{-5mm}\hspace{3cm}$$
$$\begin{array}{rl}
             & {\displaystyle + \sum_{\substack{\alpha_1,\alpha_2,\alpha_3,\alpha_4 \\ \beta_1,\beta_2,\beta_3,\beta_4}}
\frac{\omega((\alpha_1,\alpha_2,\alpha_3,\alpha_4),(\alpha_1+1,\alpha_2+1,\alpha_3+1,\alpha_4))}{\omega((\beta_1,\beta_2,\beta_3,\beta_4),
(\beta_1+1,\beta_2+1,\beta_3+1,\beta_4))}} \vspace{3pt} \\
&  \hfill e_{(\alpha_1,\alpha_2,\alpha_3,\alpha_4)} \tau e_{(\beta_1,\beta_2,\beta_3,\beta_4)}\ot e_{(\alpha_1+1,\alpha_2+1,\alpha_3+1,\alpha_4)} \tau e_{(\beta_1+1,\beta_2+1,\beta_3+1,\beta_4)} \vspace{10pt} \\
& {\displaystyle +\sum_{\substack{\alpha_1,\alpha_2,\alpha_3,\alpha_4 \\ \beta_1,\beta_2,\beta_3,\beta_4}}
\frac{\omega((\alpha_1,\alpha_2,\alpha_3,\alpha_4),(\alpha_1,\alpha_2,\alpha_3+1,\alpha_4+1))}{\omega((\beta_1,\beta_2,\beta_3,\beta_4),(\beta_1,\beta_2,\beta_3+1,\beta_4+1))}} \vspace{3pt} \\
& \hfill e_{(\alpha_1,\alpha_2,\alpha_3,\alpha_4)} \tau e_{(\beta_1,\beta_2,\beta_3,\beta_4)}\ot e_{(\alpha_1,\alpha_2,\alpha_3+1,\alpha_4+1)} \tau e_{(\beta_1,\beta_2,\beta_3+1,\beta_4+1)} \vspace{10pt} \\
 & {\displaystyle +\sum_{\substack{\alpha_1,\alpha_2,\alpha_3,\alpha_4 \\ \beta_1,\beta_2,\beta_3,\beta_4}}
\frac{\omega((\alpha_1,\alpha_2,\alpha_3,\alpha_4),(\alpha_1+1,\alpha_2+1,\alpha_3,\alpha_4+1))}{\omega((\beta_1,\beta_2,\beta_3,\beta_4),
(\beta_1+1,\beta_2+1,\beta_3,\beta_4+1))}} \vspace{-10pt} \\
& \hfill {\displaystyle e_{(\alpha_1,\alpha_2,\alpha_3,\alpha_4)} \tau e_{(\beta_1,\beta_2,\beta_3,\beta_4)}\ot   e_{(\alpha_1+1,\alpha_2+1,\alpha_3,\alpha_4+1)} \tau e_{(\beta_1+1,\beta_2+1,\beta_3,\beta_4+1)}\scalebox{1.2}{$\Bigg)$}.}
\end{array}$$

Under the identification of $\widehat{M}$ with $\widehat{P} \times \widehat{Q}$, every idempotent $e_{\phi} \in K\hspace{-1pt}M$ can\vspace{-0.5pt} be written in a unique way as $e_{\phi_1}e_{\phi_2}$ with $\phi_1 \in \widehat{P}$ and $\phi_2 \in \widehat{Q}$ such\vspace{0.5pt} that $\phi=(\phi_1,\phi_2)$. Then $(Id_{K\hspace{-1pt}S_4} \otimes \varepsilon_Q)(e_{\phi})=e_{\phi_1}\delta_{\varepsilon,\phi_2}$. In our notation with tuples this reads as:
$$\big(Id_{K\hspace{-1pt}S_4} \otimes \varepsilon_Q\big)\big(e_{(\alpha_1,\alpha_2,\alpha_3,\alpha_4)}\big)=e_{(\alpha_1,\alpha_2)}\delta_{0,\alpha_3}\delta_{0,\alpha_4}.$$
We apply $Id_A \otimes Id_{K\hspace{-1pt}S_4} \otimes \varepsilon_Q$ to the previous computation of $\Delta_J(\Psi)$ and regard the resulting element in $A \otimes K\hspace{-1pt}S_4$:
$$\begin{array}{l}
(Id_A \otimes Id_{K\hspace{-1pt}S_4} \otimes \varepsilon_Q)\Delta_J(\Psi) \vspace{10pt} \\
\hspace{3cm} \hspace{-2mm}{\displaystyle =4\scalebox{1.2}{$\Big($}\sum_{\alpha_1,\alpha_2,\beta_1,\beta_2}
\frac{\omega((\alpha_1,\alpha_2,0,0),(\alpha_1,\alpha_2,0,0))}{\omega((\beta_1,\beta_2,0,0),
(\beta_1,\beta_2,0,0))}} \vspace{2pt} \\
\hspace{3cm} \hspace{1.2cm} e_{(\alpha_1,\alpha_2,0,0)} \tau e_{(\beta_1,\beta_2,0,0)}\ot e_{(\alpha_1,\alpha_2)} \tau e_{(\beta_1,\beta_2)} \vspace{10pt} \\
\hspace{3cm} \hspace{-2mm}{\displaystyle + \sum_{\alpha_1,\alpha_2,\beta_1,\beta_2}
\frac{\omega((\alpha_1,\alpha_2,1,0),(\alpha_1+1,\alpha_2+1,0,0))}{\omega((\beta_1,\beta_2,1,0),
(\beta_1+1,\beta_2+1,0,0))}} \vspace{2pt} \\
\hspace{3cm}  \hspace{1.2cm} e_{(\alpha_1,\alpha_2,1,0)} \tau e_{(\beta_1,\beta_2,1,0)}\ot e_{(\alpha_1+1,\alpha_2+1)} \tau e_{(\beta_1+1,\beta_2+1)} \vspace{10pt} \\
\hspace{3cm} \hspace{-2mm}{\displaystyle +\sum_{\alpha_1,\alpha_2,\beta_1,\beta_2}
\frac{\omega((\alpha_1,\alpha_2,1,1),(\alpha_1,\alpha_2,0,0))}{\omega((\beta_1,\beta_2,1,1),(\beta_1,\beta_2,0,0))}} \vspace{2pt} \\
\hspace{3cm} \hspace{1.2cm} e_{(\alpha_1,\alpha_2,1,1)} \tau e_{(\beta_1,\beta_2,1,1)}\ot e_{(\alpha_1,\alpha_2)} \tau e_{(\beta_1,\beta_2)} \vspace{10pt} \\
\hspace{3cm} \hspace{-2mm}{\displaystyle +\sum_{\alpha_1,\alpha_2,\beta_1,\beta_2}
\frac{\omega((\alpha_1,\alpha_2,0,1),(\alpha_1+1,\alpha_2+1,0,0))}{\omega((\beta_1,\beta_2,0,1),
(\beta_1+1,\beta_2+1,0,0))}} \\
\hspace{3cm} \hspace{1.2cm} {\displaystyle e_{(\alpha_1,\alpha_2,0,1)} \tau e_{(\beta_1,\beta_2,0,1)}\ot
e_{(\alpha_1+1,\alpha_2+1)} \tau e_{(\beta_1+1,\beta_2+1)}\scalebox{1.2}{$\Big)$}.}
\end{array}$$

The following argument is similar to that exhibited in the uniqueness part of the proof of Proposition \ref{caseS4} (just after Table \ref{ctabS4}). We next evaluate $Id_A \otimes \chi_4$ at this element. Since $\chi_4$ is a character, we have:
$$\chi_4(e_{(\alpha_1,\alpha_2)}\tau e_{(\beta_1,\beta_2)}) = \chi_4(\tau e_{(\beta_1,\beta_2)}e_{(\alpha_1,\alpha_2)})= \chi_4(\tau e_{(\alpha_1,\alpha_2)})\delta_{\alpha_1,\beta_1}\delta_{\alpha_2,\beta_2}.$$
We also have:
$$\chi_4(\tau) = 0, \hspace{8pt} \chi_4(\tau t_1)=1, \hspace{8pt} \chi_4(\tau t_2)=-1,\hspace{8pt} \textrm{and} \hspace{8pt} \chi_4(\tau t_1t_2)=0.$$
This implies that
$$\chi_4(\tau e_{(0,0)})=0, \hspace{8pt} \chi_4(\tau e_{(0,1)})=\frac{1}{2}, \hspace{8pt} \chi_4(\tau e_{(1,0)})=-\frac{1}{2}, \hspace{8pt} \textrm{and} \hspace{7pt} \chi_4(\tau e_{(1,1)})=0.$$

Bearing all this in mind, we compute and obtain:
$$\begin{array}{rl}
x & \hspace{-2mm}=(Id_A \ot \theta)\Delta_J(\Psi) \vspace{3pt} \\
  & \hspace{-2mm}= 2\big(e_{(0,1,0,0)} \tau e_{(0,1,0,0)} - e_{(1,0,0,0)} \tau e_{(1,0,0,0)} \vspace{3pt} \\
  &  \hspace{2.5mm}+e_{(1,0,1,0)} \tau e_{(1,0,1,0)} - e_{(0,1,1,0)} \tau e_{(0,1,1,0)} \vspace{3pt} \\
  &  \hspace{2.5mm}+e_{(0,1,1,1)} \tau e_{(0,1,1,1)} - e_{(1,0,1,1)} \tau e_{(1,0,1,1)} \vspace{3pt} \\
  &  \hspace{2.5mm}+e_{(1,0,0,1)} \tau e_{(1,0,0,1)} - e_{(0,1,0,1)} \tau e_{(0,1,0,1)}\big).
\end{array}$$
Under the identification of $A$ with $K\hspace{-1pt}S_4 \otimes K\hspace{-1pt}Q$, in the latter this element is:
$$\begin{array}{rl}
x & \hspace{-2mm}= 2\big((e_{(0,1)} \tau e_{(0,1)} - e_{(1,0)} \tau e_{(1,0)}) \otimes e_{(0,0)} \vspace{3pt} \\
  &  \hspace{2.5mm}+ (e_{(1,0)} \tau e_{(1,0)} - e_{(0,1)} \tau e_{(0,1)}) \otimes e_{(1,0)} \vspace{3pt} \\
  &  \hspace{2.5mm}+ (e_{(0,1)} \tau e_{(0,1)} - e_{(1,0)} \tau e_{(1,0)}) \otimes e_{(1,1)} \vspace{3pt} \\
  &  \hspace{2.5mm}+ (e_{(1,0)} \tau e_{(1,0)} - e_{(0,1)} \tau e_{(0,1)}) \otimes e_{(0,1)}\big) \vspace{3.5pt} \\
  &  \hspace{-2mm}= 2\big(e_{(0,1)} \tau e_{(0,1)} - e_{(1,0)} \tau e_{(1,0)}\big) \otimes \big(e_{(0,0)}-e_{(1,0)}-e_{(0,1)}+e_{(1,1)}\big).
\end{array}$$
The idempotents $e_{(i,j)}$'s in the first tensorand are in $K\langle t_1,t_2\rangle$ whereas the ones in the second tensorand are in $K\langle t_3,t_4\rangle$. Notice that the element in the left-hand side tensorand is just the element $r$ of Equation \ref{elementr}. The computation done there gives:
$$\begin{array}{rl}
x & {\displaystyle \hspace{-2mm}=\frac{1}{4}\big((13)-(1234)+(23)-(1342)-(1243)+(14)-(1432)+(24)\big) \otimes (56)(78).}
\end{array}$$
Using also the computation of $r^2$ there, we obtain the square of $x$:
$$x^2 = \frac{1}{2}\big(id-(12)(34)\big) \otimes id.$$
This element, regarded as an element of $A$ is:
\begin{equation}\label{xsquare}
x^2 = \frac{1}{2}\big(id-(12)(34)\big).
\end{equation}
By construction, $x=(Id_A \ot \theta)\Delta_J(\Psi)$. Since $\theta$ and $\Psi$ are a character and cocharacter of $A$ respectively, and $X \cap A$ is a Hopf order of $A$, Proposition \ref{character} yields $x \in X \cap A$. Hence $x^2 \in X$. \par \smallskip

We repeat this argument but viewing now $S_4$ inside $S_8$ as permutations of $\{3,4,5,6\}$. Concretely, we reassign as follows: $ 1 \mapsto 3; 2 \mapsto 4; 3 \mapsto 5; 4 \mapsto 6; 5 \mapsto 1; 6 \mapsto 2;$ and $7$ and $8$ remain unchanged. Then, we obtain that
$$\frac{1}{2}\big(id-(34)(56)\big)$$
belongs to $X$. Analogously,
$$\frac{1}{2}\big(id-(56)(78)\big)$$
belongs to $X$ as well. Subtracting each of them to $id$ we get three similar elements in which the minus sign is replaced by the plus one. These elements belong to $X$ too. \par \smallskip

All idempotents $e_{(\alpha_1,\alpha_2,\alpha_3,\alpha_4)} + e_{(\alpha_1+1,\alpha_2+1,\alpha_3+1,\alpha_4+1)}$ lie in $X$ because they can be obtained from the previous elements as follows:
$$\begin{array}{l}
e_{(\alpha_1,\alpha_2,\alpha_3,\alpha_4)} + e_{(\alpha_1+1,\alpha_2+1,\alpha_3+1,\alpha_4+1)}  \vspace{5pt} \\
\hspace{1.7cm} {\displaystyle =\frac{1}{2}\big(id+(-1)^{\alpha_1+\alpha_2}t_1t_2\big)\frac{1}{2}\big(id+(-1)^{\alpha_2+\alpha_3}t_2t_3\big)\frac{1}{2}\big(id+(-1)^{\alpha_3+\alpha_4}t_3t_4
\big).}
\end{array}$$
\medskip

{\it Step 3. Replacing the twist.} Consider the following $2$-cocycle $\kappa$ on $\widehat{M}$:
$$\kappa(\varphi_1^{\alpha_1}\ldots \varphi_4^{\alpha_4},\varphi_1^{\beta_1}\ldots \varphi_4^{\beta_4}) = (-1)^{\alpha_1\beta_3+\alpha_3\beta_1}\omega(\varphi_1^{\alpha_1}\ldots \varphi_4^{\alpha_4},\varphi_1^{\beta_1}\ldots \varphi_4^{\beta_4}).$$
Observe that $\omega$ and $\kappa$ are cohomologous: defining $q:\widehat{M} \rightarrow K^{\times}$ as
$$q(\varphi_1^{\alpha_1}\ldots \varphi_4^{\alpha_4})=(-1)^{\alpha_1\alpha_3},$$
we have $\kappa=\omega \partial(q)$. \par \smallskip

Let $T$ be the twist for $K\hspace{-1pt}M$ afforded by $\kappa$. We have a Hopf algebra isomorphism \linebreak $f:(K\hspace{-1pt}S_8)_J \rightarrow (K\hspace{-1pt}S_8)_T$. We next express $T$ in the basis $\{(m,e_{\phi})\}_{m \in M, \phi \in \widehat{M}}$ of \linebreak $K\hspace{-1pt} M \otimes K\hspace{-1pt} M$. We first write it as follows:
$$T=\sum_{\phi,\psi \in \widehat{M}} \kappa(\phi,\psi)e_{\phi} \otimes e_{\psi} = \sum_{\psi \in \widehat{M}} \Big(\sum_{\phi \in \widehat{M}} \kappa(\phi,\psi)e_{\phi}\Big) \otimes e_{\psi} = \sum_{\psi \in \widehat{M}} \ell(\psi) \otimes e_{\psi},$$
where
$$\ell(\psi)=\sum_{\phi \in \widehat{M}} \kappa(\phi,\psi)e_{\phi}.$$
Evaluate $\eta \in \widehat{M}$ at $\ell(\psi)$:
$$\eta(\ell(\psi))=\sum_{\phi \in \widehat{M}} \kappa(\phi,\psi)\eta(e_{\phi})= \kappa\Big(\sum_{\phi \in \widehat{M}} \phi\eta(e_{\phi}),\psi\Big)=\kappa(\eta,\psi).$$
In the second equality we have set $\kappa$ as well for the bilinear extension of $\kappa$ to $K\hspace{-1pt}\widehat{M} \times K\hspace{-1pt}\widehat{M}$. Identify the character group of $\widehat{M}$ with $M$. Thus we get a map \linebreak $\ell: \widehat{M} \rightarrow M$. Using that $\kappa$ is a bicharacter and $\{e_{\phi}\}_{\phi \in \widehat{M}}$ a set of orthogonal idempotents, it follows that $\ell$ is indeed a group morphism. A direct computation shows that:
$$\ell(\varphi_1)=t_3, \hspace{8pt} \ell(\varphi_2)=t_1, \hspace{8pt} \ell(\varphi_3)=t_2, \hspace{8pt} \textrm{and} \hspace{8pt} \ell(\varphi_4)=t_1t_2t_3.$$
Then:\vspace{3pt}
$$\begin{array}{rl}
T & \hspace{-5pt}=  id \otimes \big(e_{(0,0,0,0)}+ e_{(1,1,1,1)}\big) + t_1 \otimes \big(e_{(0,1,0,0)}+e_{(1,0,1,1)}\big) \vspace{5pt} \\
       & \phantom{=} + t_2 \otimes \big(e_{(0,0,1,0)}+e_{(1,1,0,1)}\big) + t_3 \otimes \big(e_{(0,1,1,1)}+e_{(1,0,0,0)}\big) \vspace{5pt} \\
       & \phantom{=} + t_1t_2 \otimes \big(e_{(0,1,1,0)}+e_{(1,0,0,1)}\big) + t_1t_3 \otimes \big(e_{(0,0,1,1)}+e_{(1,1,0,0)}\big) \vspace{5pt} \\
       & \phantom{=} + t_2t_3 \otimes \big(e_{(0,1,0,1)}+e_{(1,0,1,0)}\big) + t_1t_2t_3 \otimes \big(e_{(0,0,0,1)}+e_{(1,1,1,0)}\big).
\end{array}\vspace{5pt}$$

Now, $f(X)$ is a Hopf order of $(K\hspace{-1pt}S_8)_T$. Since $f \vert_{K\hspace{-1pt}M}$ is the identity, the idempotents
$$e_{(\alpha_1,\alpha_2,\alpha_3,\alpha_4)} + e_{(\alpha_1+1,\alpha_2+1,\alpha_3+1,\alpha_4+1)}$$
belong to $f(X)$. The transpositions $t_1,t_2,$ and $t_3$ are group-like elements of $(K\hspace{-1pt}S_8)_T$. They are in $f(X)$ as well thanks to Proposition \ref{character}. Hence, $T \in f(X) \otimes_R f(X)$. Since $T=T^{-1}$, the twisted order $f(X)_{T^{-1}}$ is a Hopf order of $K\hspace{-1pt}S_8$ in virtue of Proposition \ref{twistorders}. Set $Y=f(X)_{T^{-1}}$. Again by Proposition \ref{character}, $S_8$ is contained in $Y$. We showed in \eqref{xsquare} that
$$\frac{1}{2}\big(id-(12)(34)\big)$$
belongs to $X \cap K\hspace{-1pt} M$, and hence to $Y$. Thus, the element
$$\frac{1}{2}\big(id-(12)(34)\big)+(12)(34)(45)\Big(\frac{1}{2}\big(id-(12)(34)\big)\Big)(45)=\frac{1}{2}\big(id-(354)\big)$$
belongs to $Y$. Let $\sigma=(143)(52)$. By the same reasons,
$$y:=\sigma\Big(\frac{1}{2}\big(id-(354)\big)\Big)\sigma^{-1}=\frac{1}{2}\big(id-(123)\big)$$
belongs to $Y$, and consequently, to $Y \cap (K\hspace{-1pt}S_4)$. The latter is a Hopf order of $K\hspace{-1pt}S_4$ in light of Proposition \ref{subsquo}(iii). Applying the character $\chi_4$ of Table \ref{ctabS4} and using Proposition \ref{character} for the last time we get:
$$\chi_4(y)=\frac{3}{2} \in R;$$
which implies $\frac{1}{2} \in R$.
\epf

Similarly to the proof of \cite[Corollary 2.4]{CM} one can show:

\begin{corollary}\label{complexif}
For $n \geq 4$ even the complex semisimple Hopf algebra $(\Co S_{2n})_J$ does not admit a Hopf order over any number ring.
\end{corollary}

\section{Concluding remarks}
\setcounter{equation}{0}

Nikshych showed in \cite[Corollary 4.3]{N} that if $G$ is a finite simple group, then $(\Co G)_{\Omega}$ is a simple Hopf algebra for any twist $\Omega$. Therefore, $(\Co A_n)_J$ is simple for $n \geq 5$. These were the first non-trivial examples of simple and semisimple Hopf algebras. On the other hand, Galindo and Natale proved in \cite[Theorem 3.7 and Remark 3.8]{GN} that $(\Co S_{2n})_J$ is simple for $n \geq 4$ even and not simple for $n=2$. Our results here and in \cite{CM} suggest the following questions:

\begin{question}
Let $G$ be a finite group and $\Omega$ a non-trivial twist for $\Co G$, arising from an abelian subgroup, such that $(\Co G)_{\Omega}$ is simple. Can $(\Co G)_{\Omega}$ admit a Hopf order over a number ring?
\end{question}

More generally:

\begin{question}
Let $H$ be a non-trivial simple and semisimple Hopf algebra over $\Co$. Can $H$ admit a Hopf order over a number ring?
\end{question}

Group algebras and their duals admit Hopf orders over number rings. Our results show that the twisting operation does not preserve the existence of such orders. Another class of Hopf algebras constructed from these two examples is that of lower-semisolvable Hopf algebras, defined by Montgomery and Witherspoon in \cite{MW}. Let $H$ be a Hopf algebra over $\Co$. Recall that a lower normal series for $H$ is a series of proper Hopf subalgebras
$$\Co=H_{r+1} \subset H_r \subset \ldots \subset H_2 \subset H_1=H,$$
where $H_{i+1}$ is normal in $H_i$ for each $i$. It is said that $H$ is lower-semisolvable if there is a lower normal series for $H$ in which each factor $H_i/H_iH_{i+1}^+$ is either a group algebra or the dual of a group algebra. It is natural to ask if the property of admitting a Hopf order over a number ring is inherited from the factors.

\begin{question}
Let $H$ be a semisimple Hopf algebra over $\Co$ that is lower-semisolvable. Does $H$ admit a Hopf order over a number ring?
\end{question}

Each term $H_i$ is an extension of $H_{i+1}$ and $H_i/H_iH_{i+1}^+$. We can ask more generally:

\begin{question}
Let $\Co \rightarrow A \rightarrow H \rightarrow B \rightarrow \Co$ be a short exact sequence of semisimple Hopf algebras. Suppose that $A$ and $B$ admit Hopf orders over number rings. Does $H$ admit a Hopf order as well?
\end{question}
\medskip

\subsection*{Acknowledgements}
Juan Cuadra was supported by grant MTM2017-86987-P from MICINN and FEDER and by the research group FQM0211 from Junta de Andaluc\'{\i}a. Ehud Meir was supported by the RTG 1670 ``Mathematics inspired by String theory and Quantum Field Theory''.\par \smallskip

The authors would like to thank Sonia Natale for her comments on a first version of this paper and, specially, for drawing their attention to \cite{EG1}. The authors are indebted to the attentive referee for his careful revision and apt comments, which helped to improve the original manuscript.

\end{document}